\documentclass[11pt]{article}
\usepackage{amsthm,amsmath,latexsym,amssymb}

\newcommand{\bP}{{\rm |\kern-.15em P}}
\newcommand{\Q}{\kern.3em\rule{.07em}{.65em}\kern-.3em{\rm Q}}
\newcommand{\R}{{\rm I\kern-.15em R}}
\newcommand{\D}{{\rm |\kern-.15em D}}
\newcommand{\h}{{\rm |\kern-.15em H}}
\newcommand{\C}{\kern.3em\rule{.07em}{.65em}\kern-.3em{\rm C}}
\newcommand{\T}{{\rm T\kern-.35em T}}

\theoremstyle{plain}
\newtheorem{theorem}{Theorem}[section]
\newtheorem{lemma}[theorem]{Lemma}

\newtheorem{proposition}[theorem]{Proposition}
\newtheorem{corollary}[theorem]{Corollary}

\theoremstyle{definition}
\newtheorem{definition}[theorem]{Definition}
\newtheorem{example}[theorem]{Example}

\theoremstyle{remark}
\newtheorem{remark}[theorem]{Remark}

\begin{document}
\title{Applications of Steiner symmetrization to some extremal problems in geometric function theory}
\author{Ronen Peretz}
 
\maketitle

\begin{abstract}
In this paper we investigate properties of the Steiner symmetrization in the complex plane. We use two
recursive dynamic processes in order to derive some sharp inequalities on analytic functions in the unit
disk. We answer a question that was asked by Albert Baernstein II, regarding the coefficients of circular
symmetrization. We mostly deal with the Steiner symmetrization $G$ of an analytic function $f$ in the unit
disk $U$. We pose few problems we can not solve. An intriguing one is that of the inequality
$$
\int_{0}^{2\pi} |f(re^{i\theta})|^{p}d\theta\le\int_{0}^{2\pi} |G(re^{i\theta})|^{p}d\theta,\,\,0<p<\infty
$$
which is true for $p=2$ (we prove) but can not be true for too large $p$. What is the largest such exponent
or its supremum?
\end{abstract}

\section{Some extremal problems}

\begin{definition}
Let $2\le p\le\infty$, $0<\alpha<\infty$. We define:
$$
S(p,\alpha)=\{f\,|\,f\in H(U), f\,{\rm is}\,\,{\rm univalent}\,\,{\rm in}\,U, f(0)=0, 1\le |f'(0)|, \alpha\le ||f ||_p\}
$$
where $H(U)$ is the space of all the functions that are holomorphic in the unit disc $U=\{z\in\mathbb{C}\,|\,|z|<1\}$, and
$$
||f||_{p}^{p}=\lim_{r\rightarrow 1^{-}}\frac{1}{2\pi}\int_{0}^{2\pi}|f(re^{i\theta})|^{p}d\theta.
$$
\end{definition}

\begin{definition}
Let $2\le p\le\infty$, $0<\alpha<\infty$. We define:
$$
N(p,\alpha)=\inf_{f\in S(p,\alpha)}\int_{0}^{2\pi}|f'(e^{i\theta})|d\theta.
$$
\end{definition}

\begin{proposition}
If $2\le p\le\infty$, $0<\alpha<\infty$, then there exists a function $f\in S(p,\alpha)$
such that:
$$
N(p,\alpha)=\int_{0}^{2\pi}|f'(e^{i\theta})|d\theta.
$$
\end{proposition}
\noindent
{\bf Proof.} \\
Since $g(z)=(\alpha+1)z\in S(p,\alpha)$ and $\int_{0}^{2\pi}|g'(e^{i\theta})|d\theta=2\pi(\alpha+1)$,
it follows that $N(p,\alpha)\le 2\pi (\alpha+1)$. So it will suffice to consider the following subfamily
$B(p,\alpha)=\{f\,|\,f\in S(p,\alpha),\,||f||_{\infty}\le \pi(\alpha+1)\}$ of $S(p,\alpha)$. The subfamily
$B(p,\alpha)$ of $S(p,\alpha)$ is a normal family (because it is uniformly bounded). Moreover, $B(p,\alpha)$
is a compact family. For if $f_n\in B(p,\alpha)$ and $f_n\rightarrow f$ uniformly on compact subsets of $U$
then $f\in B(p,\alpha)$ or $f\equiv 0$ ($f(0)=0$). But the condition $1\le |f'(0)|$ prohibits the second
alternative. $\qed $ \\
\begin{remark}
If in the definition of $S(p,\alpha)$ the condition $1\le |f'(0)|$ would have been dropped out, then the
claim in Proposition 1.3 would have been false. Here is an:
\begin{example}
For $p=\infty$ it is clear that $2\alpha\le N(p,\alpha)$ because for every compatible function
$f$ we have $0\in f(U)$ and there is a point $\omega\in \overline{f(U)}$ that satisfies $|\omega|=\alpha$
(because $\alpha\le ||f||_{\infty}$). Now, consider a sequence of simply connected domains $\Omega_{n}$
that satisfy $0,\alpha\in \Omega_{n}$, such that these domains become narrow as $n\rightarrow\infty$ and converge
to the closed interval on the $X$-axis, $[0,\alpha]$, and have smooth boundaries $\partial\Omega_{n}$.
By the Riemann Mapping Theorem for each $n$ there exists a conformal and onto mapping $f_{n}\,:\,U\rightarrow\Omega_{n}$
such that $f_{n}(0)=0$. Clearly we have $\lim \int_{0}^{2\pi}|f'_{n}(e^{i\theta})|d\theta=2\alpha$ and hence
if the condition $1\le |f'(0)|$ would have been dropped out from the definition of $S(p,\alpha)$
we would have had $N(\infty,\alpha)=2\alpha$, but there were no extremal function. Note that in our
construction $f_{n}\rightarrow 0$ and clearly the function $0$ is not a compatible function.
\end{example}
\end{remark}

\begin{proposition}
If $2\le p\le\infty$, $0<\alpha<\infty$ and if $f\in S(p,\alpha)$ was extremal for $N(p,\alpha)$,
then the simply connected domain $f(U)$ can have no slits.
\end{proposition}
\noindent
{\bf Proof.} \\
Let us assume to the contrary that $f(U)$ had a slit $\Gamma$. Let $D=f(U)\cup\Gamma$, i.e. $D$
is the simply connected domain we obtain from $f(U)$ by erasing the slit $\Gamma$. By the
Riemann Mapping Theorem there exits a conformal mapping $F(z)$ defined on $U$ such that
$F(U)=D$ and $F(0)=0$. We define the standard mapping $\phi\,:\,U\rightarrow U$, by the
formula, $\phi(z)=F^{-1}(f(z))$. Then $|\phi(z)|\le |z|$ and $\forall\,z\in U$ $f(z)=F(\phi(z))$.
Thus $f\prec F$, i.e. $f$ is subordinated to $F$. This implies that the following three conditions
hold true: \\
1. $F(0)=0$, $F$ conformal. \\
2. $1\le |f'(0)|\le |F'(0)|$, by the Schwarz Lemma. \\
3. $\alpha\le ||f||_{p}\le ||F||_{p}$, by Littlewood Subordination Theorem, \cite{l} or page 422 in \cite{h}. \\
These imply that $F\in S(p,\alpha)$. But clearly we have $\int_{0}^{2\pi}|F'(e^{i\theta})|d\theta<
\int_{0}^{2\pi}|f'(e^{i\theta})|d\theta$, a contradiction to the fact that $f$ is extremal for $N(p,\alpha)$. $\qed $ \\
\\
We can strengthen the Proposition 1.6:

\begin{proposition}
If $2\le p\le\infty$, $0<\alpha<\infty$ and if $f\in S(p,\alpha)$ was extremal for $N(p,\alpha)$,
then the simply connected domain $f(U)$ is a convex domain.
\end{proposition}
\noindent
{\bf Proof.} \\
Let us assume to the contrary that $f(U)$ is not a convex domain. By Proposition 1.6 it follows
that there are points $\omega_{1},\omega_{2}\in\partial f(U)$ such that $\omega_{1}\ne\omega_{2}$
and such that the open non-degenerated segment $\Gamma$ between $\omega_{1}$ and $\omega_{2}$ lies
in $\mathbb{C}-\overline{f(U)}$. Let $D$ be the simply connected domain we get by the union of $f(U)$
and the bounded domain whose boundary is the segment $[\omega_{1},\omega_{2}]$ and the corresponding 
part of $\partial f(U)$ between $\omega_{1}$ and $\omega_{2}$. From this point the proof proceeds as
that of Proposition 1.6. Namely, by the Riemann Mapping Theorem there exists a conformal mapping $F$
defined on $U$ so that $F(U)=D$ and $F(0)=0$. Then $f\prec F$ and so we have the same three conditions: \\
1. $F(0)=0$, $F$ conformal. \\
2. $1\le |f'(0)|\le |F'(0)|$, by the Schwarz Lemma. \\
3. $\alpha\le ||f||_{p}\le ||F||_{p}$, by Littlewood Subordination Theorem, \cite{l} or page 422 in \cite{h}. \\
These imply that $F\in S(p,\alpha)$. But clearly we have $\int_{0}^{2\pi}|F'(e^{i\theta})|d\theta<
\int_{0}^{2\pi}|f'(e^{i\theta})|d\theta$, a contradiction to the fact that $f$ is extremal for $N(p,\alpha)$. $\qed $ \\

\begin{proposition}
If $2\le p\le\infty$, $0<\alpha<\infty$ and if $f\in S(p,\alpha)$ was extremal for $N(p,\alpha)$,
and if $|a|<1$ then either:
$$
|f'(a)|\le\frac{1}{1-|a|^2},
$$
or
$$
\int_{0}^{2\pi}|f(e^{i\theta})-f(a)|^{p}\left(\frac{1-|a|^{2}}{|e^{i\theta}-a|^{2}}\right)d\theta\le\alpha^{p}.
$$
In particular for $a=0$: either $|f'(0)|=1$ or $\int_{0}^{2\pi}|f(e^{i\theta})|^{p}d\theta=||f||_{p}^{p}=\alpha^{p}$.
\end{proposition}
\noindent
{\bf Proof.} \\
If
$$
\phi(z)=\frac{z+a}{1+\overline{a}z}
$$
then
$$
\phi'(z)=\frac{1-|a|^2}{(1+\overline{a}z)^2}.
$$
We define $F(z)=f(\phi(z))-f(a)$. Then we have: \\
1. $\int_{0}^{2\pi}|F'(e^{i\theta})|d\theta=\int_{0}^{2\pi}|f'(e^{i\theta})|d\theta$ because the images $F(U)$
and $f(U)$ are congruent. \\
2. $F'(0)=(1-|a|^2)f'(a)$. \\
\\
We also have the identity:
$$
\int_{0}^{2\pi}|F(e^{i\psi})|^{p}d\psi=\int_{0}^{2\pi}\left|f\left(\frac{e^{i\psi}+a}{1+\overline{a}e^{i\psi}}\right)-f(a)\right|^{p}d\psi.
$$
We make a change of the integration variable:
$$
e^{i\theta}=\frac{e^{i\psi}+a}{1+\overline{a}e^{i\psi}},\,\,\,\,\,d\psi=\frac{1-|a|^2}{|e^{i\theta}-a|^2}d\theta,
$$
and we rewrite the above identity as follows: \\
\\
3.
$$
\int_{0}^{2\pi}|F(e^{i\theta})|^{p}d\theta=\int_{0}^{2\pi}|f(e^{i\theta})-f(a)|^{p}\left(\frac{1-|a|^2}{|e^{i\theta}-a|^{2}}\right)d\theta.
$$
Now 1, 2 and 3 above imply that: \\
\\
4. If we have both: $1\le (1-|a|^{2})|f'(a)|$ and
$$
\alpha^{p}\le\int_{0}^{2\pi}|f(e^{i\theta})-f(a)|^{p}\left(\frac{1-|a|^{2}}{|e^{i\theta}-a|^{2}}\right)d\theta,
$$
then $F\in S(p,\alpha)$ and hence $F$ is extremal for $N(p,\alpha)$. If both the inequalities in 4 are sharp (none of
them is an equality), then there is an $M>1$ such that $M^{-1}F(z)\in S(p,\alpha)$, which contradicts the fact that
$F(z)$ is extremal for $N(p,\alpha)$. Thus at least one of the two inequalities in 4 is in fact an equality and the
proposition follows. $\qed $

\section{Facts about symmetrizations}
Motivated by the desire to solve the family of the extremal problems $N(p,\alpha)$, we will discuss in this
section properties of symmetrizations of functions. Specifically we will consider symmetrizations that were
introduced by P\"olya and by Steiner. We will recall results from the paper \cite{b}.

\begin{definition}
Let $D$ be a domain in the Riemann sphere $\mathbb{C}\cup\{\infty\}$. The circular symmetrization of $D$ is the 
domain $D^{*}$ that is defined as follows: for each $t\in (0,\infty)$ we define $D(t)=\{\theta\in [0,2\pi]\,|\,te^{i\theta}\in D\}$.
If $D(t)=[0,2\pi]$ then the intersection of $D^{*}$ with the circle $|z|=t$ is the full circle. If $D(t)=\emptyset$
then the intersection of $D^{*}$ with the circle $|z|=t$ is the empty set $\emptyset$. If $D(t)$ is a non trivial subset of
$[0,2\pi]$ which has the measure $|D(t)|=\alpha'$, then the intersection of $D^{*}$ with the circle $|z|=t$ is the unique
circular arc given by $\{te^{i\theta}\,|\,|\theta|< \alpha'/2\}$. Finally $D^{*}$ contain the point $0$ ($\infty$)
if and only if $D$ contains the point $0$ ($\infty$).
\end{definition}
\noindent
Section (j) of the paper \cite{b} includes a proof of an important principle in symmetrization: \\
Let $f\in H(U)$ and let us denote $D=f(U)$. Let $D_0$ be a simply connected domain that contains $D^{*}$, and let us
assume that $D_0$ is not the full complex plane ($\mathbb{C}$). Let $F$ be a conformal mapping of $U$ onto $D_0$
that satisfies $F(0)=|f(0)|$. The following result is proved in \cite{b}: \\
\\
{\bf Theorem 6. (\cite{b})} {\it
If  $\Phi$ is a convex non-decreasing function on $(-\infty,\infty)$, $f\in H(U)$ and $F$ as above, then for all
$0\le r<1$ we have:
$$
\int_{-\pi}^{\pi}\Phi(\log |f(re^{i\theta})|)d\theta\le\int_{-\pi}^{\pi}\Phi(\log |F(re^{i\theta})|)d\theta.
$$
}
If we choose in Theorem 6 above, $\Phi(x)=e^{2x}$ and assume that we have the following expansions: 
$f(z)=\sum_{n=0}^{\infty}a_{n}z^{n}$ and $F(z)=\sum_{n=0}^{\infty}A_{n}z^{n}$, then we obtain the inequality
$\sum_{n=0}^{\infty}|a_{n}|^{2}r^{2n}\le\sum_{n=0}^{\infty}|A_{n}|^{2}r^{2n}$ for $0\le r<1$. By the definition
of $F$ we have $|A_{0}|=|a_{0}|$, thus if we subtract $|A_{0}|^{2}$ from both sides of the inequality and divide by $r^{2}$
and than take $r\rightarrow 0^{+}$ we obtain $|f'(0)|\le |F'(0)|$, a classical result of Walter Hayman. If $f$
is one-to-one in $U$ then both $D$ and $D^{*}$ are simply connected and we can take $F$ to be a conformal mapping
from $U$ onto $D^{*}$ for which $F(0)=|f(0)|$.

At the end of section (k) in \cite{b} the author asks if the following is true for all $n$: $|a_{n}|\le |A_{n}|$?
Is the following weaker set of inequalities true: $\sum_{k=0}^{n}|a_{k}|^{2}\le\sum_{k=0}^{\infty}|A_{k}|^{2}$?
A. Baernstein II, remarks that these last inequalities if true, would prove a conjecture of Littlewood: If $f$
is one-to-one and analytic in $U$ and if $f(z)\ne 0$, for $z\in U$, then for each $n>1$ we have: $a_{n}\le 4n|a_{0}|$. \\
(We give a proof for this assertion, for the reader's convenience. Assume $a_{0}=1$. We define $g(z)=\sqrt{f(z)}$ and
$h(z)=-g(z)$. Then $g,h\in H(U)$ (because $f(z)$ never vanishes in $U$). Both $g$ and $h$ are one-to-one in $U$ (because
$g(z_{1})=g(z_{2})\Rightarrow f(z_{1})=f(z_{2})\Rightarrow z_{1}=z_{2}$). Let us suppose that 
$g(z)=\sum_{n=0}^{\infty}\alpha_{n}z^{n}$, where $\alpha_{0}=1$. We denote $D=g(U)$. If $\xi_{0}\in D$ then $-\xi_{0}\not\in D$
(because $g(U)\cap h(U)=\emptyset$). Hence for any $0<t<\infty$ we must have $|D(t)|\le\pi$. We conclude that
for the symmetrization we have: $D^{*}\subseteq\{z\in\mathbb{C}\,|\,\Re\{z\}\ge 0\}$. Let $G(z)$ be a conformal mapping
from $U$ onto $D^{*}$ such that $G(0)=g(0)=1$. Then for all $|z|<1$ we have $\Re\{z\}\ge 0$ and hence if
$G(z)=1+\sum_{n=1}^{\infty}B_{n}z^{n}$ then by a result of Constantine Carath\'eodory we have $|B_{n}|\le 2$ for all
$n=1,2,3,\ldots $. We conclude that $\sum_{k=1}^{n}|B_{k}|^{2}\le 4n$ which implies (assuming Baernstein's assumption above)
that $\sum_{k=0}^{n}|\alpha_{k}|^{2}\le 4n$. We recall that $f=g^{2}$ and hence:
$$
|a_{n}|^{2}=|\alpha_{0}\alpha_{n}+\alpha_{1}\alpha_{n-1}+\ldots+\alpha_{n}\alpha_{0}|^{2}\le
$$
$$
\le(2|\alpha_{n}|+|\alpha_{1}||\alpha_{n-1}|+\ldots+|\alpha_{n-1}||\alpha_{0}|)^{2}\le (4+4(n-1))\cdot 4n=(4n)^{2}.
$$
$\qed $)

Concerning the first question posed by Albert Baernstein II (above): Let $f$ be a conformal function defined on $U$. 
We assume that $f(U)$ has a finite area. Let us denote $D=f(U)$, and let $F$ be a conformal mapping of $U$ onto the 
symmetrization $D^{*}$ such that $F(0)=|f(0)|$. Let us denote by $S(D)$ and by $S(D^{*})$ the areas of the respective 
domains. We will use $td\phi\cdot dt$ for the area element in polar coordinates. Then we have the identities:
$$
S(D)=\int_{0}^{\infty}\int_{D(t)}td\phi\cdot dt=\int_{0}^{\infty}t|D(t)|dt,\,\,\,\,S(D^{*})=\int_{0}^{\infty}t|D^{*}(t)|dt.
$$
By the definition of $D^{*}$ it follows that for all $0\le t<\infty$ we have $D(t)=D^{*}(t)$ and hence $S(D)=S(D^{*})$
(the well-known fact that circular symmetrization is an area preserving transformation). On the other hand we have
$S(D)=\int_{0}^{1}\int_{0}^{2\pi}r|f'(re^{i\theta})|^{2}d\theta dr$ and if $f(z)=\sum_{n=0}^{\infty}a_{n}z^{n}$ and
$F(z)=\sum_{n=0}^{\infty}A_{n}z^{n}$ then we obtain the well known formulas: $S(D)=\pi\sum_{n=0}^{\infty}n|a_{n}|^{2}$, and
$S(D^{*})=\pi\sum_{n=0}^{\infty}n|A_{n}|^{2}$. We conclude that $\sum_{n=0}^{\infty}n|a_{n}|^{2}=\sum_{n=0}^{\infty}
n|A_{n}|^{2}$. We recall that by the definition of $F$ we have $A_{0}=|a_{0}|$ and by Hayman's result $|a_{1}\le |A_{1}|$
and so either $|a_{n}|=|A_{n}|$ for $n=0,1,2,\ldots $ or there exist $1\le n_{1},n_{2}$ so that $|a_{n_{1}}|<|A_{n_{1}}|$,
and $|A_{n_{2}}|<|a_{n_{2}}|$. We proved the following:
\begin{theorem}
If $f(z)=\sum_{n=0}^{\infty}a_{n}z^{n}$ is analytic, one-to-one in $U$ and $f(U)$ has a finite area. If 
$F(z)=\sum_{n=0}^{\infty}A_{n}z^{n}$ is the circular symmetrization of $f(z)$, then we have: 
$\sum_{n=0}^{\infty}n|a_{n}|^{2}=\sum_{n=0}^{\infty}n|A_{n}|^{2}$ and either for all $n=0,1,2,\ldots $
we have $|a_{n}|=|A_{n}|$, or there exist $1\le n_{1},n_{2}$ such that $|a_{n_{1}}|<|A_{n_{1}}|$ and
$|A_{n_{2}}|<|a_{n_{2}}|$.
\end{theorem} 
\noindent
Theorem 2.2 answers the problem mentioned above that was raised by Albert Baernstein II. The answer in negative. \\
Next, let $f(z)=\sum_{n=0}^{\infty} a_{n}z^{n}$, $z\in U$, be analytic and one-to-one, and assume that
$\partial D=\partial f(U)$ is rectifiable.

\begin{remark}
Let $0<a<b$, $0<h$. Then the sum of the lengths of the legs of a trapezoidal with bases of lengths $a$ and $b$
and height of length $h$ is minimal, if and only if the legs are equal to one another. The same fact is true also
for circular a trapezoidal which has bases (of lengths $a$ and $b$) located on two concentric circles of radii
$R$ and $R+h$.
\end{remark}
\noindent
Let $\gamma$ be a rectifiable curve. Let the (finite) length of $\gamma $ be denoted by $l(\gamma)$. We chose a finite
number of points on $\gamma $ and join successive points by straight segments. The result is a polygonal curve
that is composed of the straight segments $\gamma_{1},\gamma_{2},\ldots,\gamma_{n}$ (we assume they come geometrically
one after the other). The length of the polygonal curve is the finite sum $l_{n}(\gamma)=\sum_{j=1}^{n}l(\gamma_{j})$.
When we refine the division points and take $n\rightarrow\infty$ we obtain $l(\gamma)=\lim_{n\rightarrow\infty}l_{n}(\gamma)$.
By the definition of $D^{*}$ and by remark 2.3  it follows that in order to compute the lengths $l(\partial D)$ and
$l(\partial D^{*})$, if at each approximation step we take our division points to be the intersection points of $\partial D$
($\partial D^{*}$) with sets of concentric circles centered at the origin and of radii $0<t_{1}<t_{2}<\ldots<t_{n}$, so that
$0<t_{j+1}-t_{j}<\epsilon_{n}$, $j=1,2,\ldots,n-1$, $\epsilon_{n}\rightarrow_{n\rightarrow\infty} 0^{+}$, we obtain the 
inequalities: $l_{n}(\partial D^{*})\le l_{n}(\partial D)$. Hence when $n\rightarrow\infty$ we obtain the well known
fact $l(\partial D^{*})\le l(\partial D)$ (circular symmetrization reduces the perimeter). On the other hand we have
the two identities:
$$
l(\partial D)=\int_{0}^{2\pi}|f'(e^{i\theta})|d\theta,\,\,\,\,l(\partial D^{*})=\int_{0}^{2\pi}|F'(e^{i\theta})|d\theta.
$$
Here $f$ is a conformal mapping defined on $U$ with the image $D=f(U)$ and $F$ is a conformal mapping defined on $U$ with
the image $F(U)=D^{*}$ (which is simply connected as well).

\begin{definition}
Let $D$ be a domain in the Riemann sphere $\mathbb{C}\cup\{\infty\}$. The Steiner symmetrization of $D$ is the 
domain $D_{S}^{*}$ that is defined as follows: for each $t\in (-\infty,\infty)$ we define $D(t)=\{t+iy\in\mathbb{C}\,|\,t+iy\in D\}$.
If $D(t)=\emptyset$ then the intersection of $D_{S}^{*}$ with the vertical line $H(t)=\{z\in\mathbb{C}\,|\,\Re\{z\}=t\}$ is the empty set
$\emptyset$. If $D(t)=D\cap H(t)$ is a non-empty and has the measure $|D(t)|=\alpha'$, then the intersection of $D_{S}^{*}$ with the line $H(t)$ is the unique segment or line given by $\{t+iy\in\mathbb{C}\,|\,|y|<\alpha'/2\}$.
\end{definition}
\noindent
Let $G$ be a conformal mapping defined on $U$ with the image $G(U)=D_{S}^{*}$ (which is simply connected). We assume as above that the 
boundary curves $\partial D=\partial f(U)$ and $\partial D_{S}^{*}=\partial G(U)$ are rectifiable. Then the remark
above shows that we also have $l_{n}(\partial D_{S}^{*})\le l_{n}(\partial D)$, where this time we take our division points to be the intersection points of $\partial D$
($\partial D_{S}^{*}$) with sets of parallel vertical lines of $x$-coordinates $-\infty<t_{1}<t_{2}<\ldots<t_{n}<\infty$, so that
$0<t_{j+1}-t_{j}<\epsilon_{n}$, $j=1,2,\ldots,n-1$, $\epsilon_{n}\rightarrow_{n\rightarrow\infty} 0^{+}$. We have the identity:
$$
l(\partial D_{S}^{*})=\int_{0}^{2\pi}|G'(e^{i\theta})|d\theta.
$$
\noindent
We just proved the well known:

\begin{theorem}
If $f\in H(U)$, $f$ is one-to-one, and the boundary curve $\{f(e^{i\theta})\,|\,0\le\theta<2\pi\}$ is rectifiable. If $F$
is the circular symmetrization (P\"olya) of $f$ and if $G$ is the Steiner symmetrization of $f$ ($F(0)=G(0)=|f(0)|$), then
the boundary curves $\{F(e^{i\theta})\,|\,0\le\theta<2\pi\}$ and $\{G(e^{i\theta})\,|\,0\le\theta<2\pi\}$ are rectifiable, and
we have the following two inequalities:
$$
\int_{0}^{2\pi}|G'(e^{i\theta})|d\theta\le \int_{0}^{2\pi}|f'(e^{i\theta})|d\theta,\,\,\,\,\,
\int_{0}^{2\pi}|F'(e^{i\theta})|d\theta\le \int_{0}^{2\pi}|f'(e^{i\theta})|d\theta.
$$
\end{theorem}

\begin{remark}
We comment the two items below in connection to the second question of Baernstein and the possible proof of the
conjecture of Littlewood on non-vanishing one-to-one analytic functions in $U$. We supply proofs to these well known
facts, for the convenience of the reader:\\
\\
1. If $F(z)$ is analytic and one-to-one in $U$ and satisfies $0<\Re\{F(z)\}$, then $(F(z))^{2}$ is one-to-one in $U$. \\
\\
{\bf Proof.} \\
$(F(z_{1}))^{2}=(F(z_{2}))^{2}\,\Rightarrow\,\,F(z_{1})=F(z_{2})\,\,{\rm or}\,\,F(z_{1})=-F(z_{2})$. But if
$F(z_{1})=-F(z_{2})$ then $\Re\{F(z_{1})\}\cdot\Re\{F(z_{2})\}< 0$ which proves that $F(z_{1})=F(z_{2})$
and hence $z_{1}=z_{2}$. $\qed $ \\
\\
2. If $f\in H(U)$ is one-to-one in $U$ and $f(z)\ne 0$ in $U$ and if $F(z)$ is the circular symmetrization of
$\sqrt{f(z)}$, then $(F(z))^{2}$ is the circular symmetrization of $f(z)$ (and vice versa). \\
\\
{\bf Proof.} \\
By remark 1 above $(F(z))^{2}$ is one-to-one which implies the claim. $\qed $
\end{remark}
\noindent
We need one more well known fact, this time on the Steiner symmetrization. Our proof
uses the fact that we can regard Steiner symmetrization as a limiting case (at infinity)
of circular (P\"olya) symmetrization.

\begin{theorem}
If $f$ is analytic and one-to-one in $U$ and if $G$ is the Steiner symmetrization of $f$, then for any $r$,
$0\le r<1$ we have the inequality:
$$
\int_{0}^{2\pi}|f(re^{i\theta})|^{2}d\theta\le \int_{0}^{2\pi}|G(re^{i\theta})|^{2}d\theta.
$$
\end{theorem}
\noindent
{\bf Proof.} \\
Let $0<M$. Let $F_{M}(z)$ be the circular symmetrization of the shifted function $M+f(z)$. Let us
denote $G_{M}(z)=F_{M}(z)-M$. We will use the theorem of Baernstein (\cite{b}, Theorem 6) cited above. 
For any $0<p<\infty$:
$$
\int_{0}^{2\pi}|M+f(re^{i\theta})|^{p}d\theta\le \int_{0}^{2\pi}|F_{M}(re^{i\theta})|^{p}d\theta.
$$
This can be re-written as follows:
$$
\int_{0}^{2\pi}|M+f(re^{i\theta})|^{p}d\theta\le \int_{0}^{2\pi}|M+G_{M}(re^{i\theta})|^{p}d\theta.
$$
which proves the following:
$$
\int_{0}^{2\pi}\left(1+\frac{p}{M}\Re\{f(re^{i\theta})\}+\frac{1}{M^{2}}\left\{\frac{p^{2}}{4}|f(re^{i\theta})|^{2}+
2\left(\begin{array}{c} p/2 \\ 2 \end{array}\right)\Re\{(f(re^{i\theta}))^{2}\}\right\}\right)d\theta\le
$$
$$
\le\int_{0}^{2\pi}\left(1+\frac{p}{M}\Re\{G_{M}(re^{i\theta})\}+\frac{1}{M^{2}}\left\{\frac{p^{2}}{4}|G_{M}(re^{i\theta})|^{2}+
2\left(\begin{array}{c} p/2 \\ 2 \end{array}\right)\Re\{(G_{M}(re^{i\theta}))^{2}\}\right\}\right)d\theta+
$$
$$
+o\left(\frac{1}{M^{2}}\right).
$$
But $\lim_{M\rightarrow\infty}G_{M}=G$ the Steiner symmetrization of $f$, uniformly on compact subsets of $U$. This proves
(on taking $M\rightarrow\infty$) that:
$$
\int_{0}^{2\pi}|f(re^{i\theta})|^{2}d\theta\le \int_{0}^{2\pi}|G(re^{i\theta})|^{2}d\theta.\,\,\,\,\,\qed
$$

\begin{remark}
The analog of the theorem of Baernstein (Theorem 6, quoted above) to Steiner symmetrization is false, i.e. the statement 
$\int_{0}^{2\pi} |f(re^{i\theta})|^{p}d\theta\le\int_{0}^{2\pi} |G(re^{i\theta})|^{p}d\theta$ for $0<p<\infty$
is false. It is easy to construct a counterexample. Thus $p=2$ is an exception. This naturally raises a question. \\
\\
{\bf A problem (the means of the Steiner symmetrization):} \\
Find the largest possible $2\le p_{0}$ for which $\forall\,0\le r<1$ and $\forall\,p\le p_{0}$ (or $\forall\,p<p_{0}$)
we have the inequality:
$$
\int_{0}^{2\pi}|f(re^{i\theta})|^{p}d\theta\le \int_{0}^{2\pi}|G(re^{i\theta})|^{p}d\theta,
$$
for any $f\in H(U)$, one-to-one in $U$, with $f(0)=0$, where $G$ is the Steiner symmetrization of $f$. \\
\\
By Theorem 2.7 above we know that $2\le p_{0}$ and that the inequality is valid for $p=2$.
\end{remark}
\noindent
We conclude this section with a simple demonstration of the type of reductions we can make so far concerning
the solution of the family of problems $N(p,\alpha)$.

\begin{proposition}
If $2\le p\le\infty$, $0<\alpha<\infty$, and if $f\in S(p,\alpha)$ is an extremal function for $N(p,\alpha)$
then we may assume that the domain $f(U)$ is circular symmetric (P\"olya symmetric).
\end{proposition}
\noindent
{\bf Proof.} \\
Let the function $F$ be the P\"olya symmetrization of the extremal function $f$. Then the following three
properties hold true: \\
1. $F\in H(U)$, $F(0)=0$, $F$ is one-to-one in $U$. \\
2. $1\le |f'(0)|\le |F'(0)|$, by a classical result of Walter Hayman. \\
3. $\alpha\le ||f||_{p}\le ||F||_{p}$, by Littlewood Subordination Theorem, \cite{l} or page 422 in \cite{h}. \\
As in Theorem 2.5, it follows that:
$$
\int_{0}^{2\pi}|F'(e^{i\theta})|d\theta\le\int_{0}^{2\pi}|f'(e^{i\theta})d\theta.
$$
But by the assumption, the function $f$ is an extremal function for the problem $N(p,\alpha)$, and so it
follows that also the function $F$ is extremal for the problem $N(p,\alpha)$. $\qed $ \\

\section{A solution of the problem $N(2,\alpha)$}
We will present a solution of the simple case $N(2,\alpha)$. Already here we will conclude a
few interesting conclusions. For example, we will be able to prove the convergence of certain infinite
products of geometric parameters which we can not easily explicitly compute. The main argument will be based
on the compactness of the family $S(2,\alpha)$. Let $f\in S(2,\alpha)$, $0<\alpha<\infty$. We will outline
a simple recursive process that constructs in each phase a function $g$ which satisfies the two properties: \\
1. $g\in S(2,\alpha)$. \\
2. $\int_{0}^{2\pi}|g'(e^{i\theta})|d\theta\le\int_{0}^{2\pi}|f'(e^{i\theta})|d\theta$. \\
\\
In fact the $1$-norm of $g'(e^{i\theta})$ will be smaller than or equals to the same norm of the function
constructed in the previous step. \\
{\bf Step 1:} We choose an angle $\phi $ and define $\overline{f}_{\phi}(z)=e^{i\phi}f(z)$. \\
{\bf Step 2:} We compute the function $g_{\phi}(z)$, the Steiner symmetrization of $\overline{f}_{\phi}$. \\
\\
We note that $g_{\phi}$ has the following 4 properties: \\
$g_{\phi}\in H(U)$, $g_{\phi}$ is one-to-one and $g_{\phi}(0)=0$. \\
$1\le |\overline{f}_{\phi}'(0)|\le |g_{\phi}'(0)|$, by the result of Walter Hayman mentioned before. \\
$\alpha\le ||\overline{f}_{\phi}||_{2}\le ||g_{\phi}||_{2}$, by Theorem 2.7. \\
$\int_{0}^{2\pi}|g_{\phi}'(e^{i\theta})|d\theta\le\int_{0}^{2\pi}|\overline{f}_{\phi}'(e^{i\theta})|d\theta$, by Theorem 2.5. \\
{\bf Step 3:} We compute the following number:
$$
c_{\phi}=\max\left(\frac{1}{|g_{\phi}'(0)|},\frac{\alpha}{||g_{\phi}||_{2}}\right),
$$
and then we have $0<c_{\phi}\le 1$. We compute $g(z)=c_{\phi}g_{\phi}(z)$.

\begin{remark}
In less formality we can describe the process that was outlined above as follows: \\
1. Rotate $f(U)$ about 0 (in an angle $\phi$). \\
2. Perform Steiner symmetrization of the rotated domain. \\
3. Shrink the domain that was obtained, by an optimal factor $c_{\phi}$, where $0<c_{\phi}\le 1$.
\end{remark}

\begin{definition}
Let $D$ be a domain that contains the origin, $0$. The $\phi$-deformation of $D$ is the domain $D_{\phi}$
which we get by the Steiner symmetrization of $e^{i\phi}D$. In other words $D_{\phi}$ is the resulting 
domain after steps 1 and 2 in the process we described above.
\end{definition}

\begin{definition}
If $f\in S(2,\alpha)$ and if $\phi$ is a real number, we will denote the function $g$ which is obtained
after executing the process above (steps 1, 2 and 3), by $g=f_{\phi}$. The number $c_{\phi}$ will be
called \underline{the shrinking factor}.
\end{definition}
\noindent
If $\{\phi_{n}\}_{n=1}^{\infty}$ is a sequence of real numbers, then \underline{the corresponding sequence}
\underline{of shrinking factors} will be defined to be the sequence of the shrinking factors we obtain by executing
the following iterative process: \\
$\{f,f_{\phi_1},(f_{\phi_1})_{\phi_2},((f_{\phi_1})_{\phi_2})_{\phi_3},\ldots\}$. \\
\\
If $\{D_{n}\}_{n=1}^{\infty}$ is a sequence of bounded domains, and if $D$ is a bounded domain, then we will say
that the sequence of the domains converges to to the domain $D$ and write $D_{n}\rightarrow D$, if for any
$\epsilon>0$ there exists a number $N$ such that for all $n>N$ we will have ${\rm distance}(\partial D_{n},
\partial D)<\epsilon$. We have the following surprising fact:

\begin{theorem}
If $0<\alpha<\infty$, and if $f\in S(2,\alpha)$ and also $\{\phi_{n}\}_{n=1}^{\infty}$ is any sequence
of real numbers, then we have the double inequality:
$$
0<\prod_{n=1}^{\infty}c_{\phi_n}\le 1.
$$
\end{theorem}
\noindent
{\bf Proof.} \\
For each $n$, $n=1,2,3,\ldots $, we have $0<c_{\phi_n}\le 1$. Hence $0\le\prod_{n=1}^{\infty} c_{\phi_n}\le 1$.
If $\prod_{n=1}^{\infty} c_{\phi_n}=0$ then $f_{\phi_n}\rightarrow 0$ uniformly on compact subsets of $U$.
This, however, contradicts the compactness of the family $S(2,\alpha)$. $\qed $ \\

\begin{theorem}
Let $D$ be a bounded domain that contains the origin, $0$. Then, there exists a disk $B$ whose center
is the origin, $0$ and there exists a sequence of $\phi_{n}$-deformations of $D$ that will be denoted by
$\{D_{n}\}_{n=1}^{\infty}$ so that $D_{n}\rightarrow B$, (The assumption means that $D_{1}$ is the
$\phi_1$-deformation of $D$, and $D_{n+1}$ is the $\phi_{n+1}$-deformation of $D_{n}$). Moreover, the
disk $B$ is unique in the sense that if $\{D_{n}'\}_{n=1}^{\infty}$ is the sequence of $\phi_{n}'$-deformations of $D$
that satisfies $D_{n}'\rightarrow B'$, then $d(B)=d(B')$
\end{theorem}
\noindent
{\bf Proof.} \\
Firstly, we note that, if $D$ is a bounded domain that contains the origin, $0$, and if $\phi $ is real
number, then the $\phi$-deformation of $D$, $D_{\phi}$, has its diameter smaller than or 
equals to the diameter $D$, i.e. $d(D_{\phi})\le d(D)$. Moreover, $d(D_{\phi})<d(D)$ always except for 
the case in which the diameter of $e^{i\phi}D$ is parallel to one of the axis of the coordinates. We denote
by $A$ the set of all the domains $\tilde{D}$ for which there exists a finite sequence of real numbers 
$\{\alpha_{n}\}_{n=1}^{N}$ such that if $D_{1}=D_{\phi_1}$, $D_{n+1}=(D_{n})_{\alpha_{n+1}}$,
$1\le n\le N-1$, then $\tilde{D}=D_{N}$. We define $d=\inf\{d(\tilde{D})\,|\,\tilde{D}\in A\}$. We claim that
$0<d<\infty$. For $0\in D$ and hence there is a disk with a positive-radius $\epsilon/2$ about the origin, $0$
is contained in $D$. Hence $\epsilon\le d \le d(D)$. Moreover, there exists a bounded domain $B$ that 
contains the origin, $0$ such that $d=d(B)$ and such that there is a sequence of $\phi_{n}$-deformations
of $D$, $\{ D_{n}\}_{n=1}^{\infty}$ that satisfy $D_{n}\rightarrow B$. We will now prove that $B$ is a
disk centered at $0$: The domain $B$ has a diameter in every direction (for a $\phi$-deformation properly shrinks
diameters except for the event in which the rotated domain has a diameter parallel to one of the axis
of the coordinates. Since $B$ has a minimal diameter, there is no way to shrink its diameter properly
using a $\phi$-deformation. The claim follows). A domain which has a diameter in every direction is
a disk (for the triangle inequalities imply that any pair of diameters must intersect through their
common mid-point). Finally, if $B$ is a disk and if $0\in B$, then using a single $\phi$-deformation
we can bring $B$ to be a disk (with the same diameter) whose center is the origin, $0$. \\
Let us prove the uniqueness claim: we note that a $\phi$-deformation (just like Steiner symmetrization)
is an area preserving transformation. Hence if $2R-d(B')$ then $\pi R^{2}=\int\int_{D}dxdy$ and we conclude
that $d(B')=d(B)=2\sqrt{\int\int_{D}dxdy/\pi}$. $\qed $ \\
\\
We shall now solve the problem $N(2,\alpha)$.

\begin{theorem}
If $0<\alpha<\infty$, and if $G(z)=\max\{1,\alpha\}z$, then $G\in S(2,\alpha)$ and we have:
$$
N(2,\alpha)=\int_{0}^{2\pi}|G'(e^{i\theta})|d\theta=2\pi\max\{1,\alpha\}.
$$
\end{theorem}
\noindent
{\bf Proof.} \\
By Theorem 3.4 and Theorem 3.5 it follows that among those functions that minimize, there is one 
that conformally maps $U$ onto a disk centered at the origin, $0$. Such a mapping has the form 
$cz$ and since it should belong to the family $S(2,\alpha)$, we necessarily have: $\max\{1,\alpha\}\le |c|$.
But by our assumption the function solves (by minimizing) the extremal problem. We conclude that
$|c|=\max\{1,\alpha\}$. $\qed $ \\

\begin{remark}
Our definition of the family $S(2,\alpha)$ requires the condition $1\le |f'(0)|$ in order to obtain
a compact family. If instead of that condition we had $\epsilon\le |f'(0)|$ for some fixed $0<\epsilon$,
we could have solved the corresponding problem $N(2,\alpha)$ similarly, except that this time
our multiplier had to be $\max\{\epsilon,\alpha\}$. In particular for small enough $\epsilon$ the
solution would have been $G(z)=\alpha z$. We conclude that if $f\in H(U)$, $f$ is one-to-one, $f(0)=0$,
and if we denote $\alpha=||f||_{2}$ then if $\alpha<\infty$ we can solve the problem $N(2,\alpha)$
with a small enough $\epsilon>0$ (meaning $\epsilon<\alpha$) and obtain a solution $\alpha z$. This 
helps in proving the following inequality:
\end{remark}

\begin{theorem}
Let $f\in H(U)$, $f(0)=0$. Then for each $0\le r<1$ we have the following inequality:
$$
\left(\frac{1}{2\pi}\int_{0}^{2\pi} |f(re^{i\theta})|^{2}d\theta\right)^{1/2}\le\frac{r}{2\pi}\int_{0}^{2\pi}
|f'(re^{i\theta})|d\theta.
$$
In particular we have $||f||_{2}\le ||f'||_{1}$. Both inequalities above are sharp.
\end{theorem}
\noindent
{\bf Proof.} \\
Let us start by assuming that we already proved the first inequality for functions $f\in H(U)$,
for which $f(0)=0$, that are also one-to-one. Let $g\in H(U)$ satisfy $g(0)=0$ and let $0<r<1$.
We consider the conformal mapping $f\in H(U)$ such that $f(0)=0$ that satisfies $f(U)=g(U)$. The
function $g$ is subordinate to the function $f$ (which means that there exists a $w\in H(U)$,
$w(0)=0$, $|w(z)|<1$ so that $g(z)=f(w(z))$). By a theorem of Littlewood, \cite{l} or page 422 in \cite{h}, we have the following 
inequality:
\begin{equation}
\label{eq1}
\frac{1}{2\pi}\int_{0}^{2\pi}|g(re^{i\theta})|^{2}d\theta\le\frac{1}{2\pi}\int_{0}^{2\pi}|f(re^{i\theta})|^{2}d\theta.
\end{equation}
Since we assumed that $f$ is conformal, it follows by our initial assumption (at the beginning of the 
proof) that the following inequality is true:
\begin{equation}
\label{eq2}
\left(\frac{1}{2\pi}\int_{0}^{2\pi} |f(re^{i\theta})|^{2}d\theta\right)^{1/2}\le\frac{r}{2\pi}\int_{0}^{2\pi}
|f'(re^{i\theta})|d\theta.
\end{equation}
Finally (by the extension theorem of the Riemann Theorem) the holomorphic mapping maps the boundary of the 
domain onto the boundary of range domain (the boundaries are smooth enough and we can use
a theorem of Constantine Carath\'eodory). Since $f$ is (by assumption) a conformal mapping it traces the
boundary $\partial g(rU)$ once, while the mapping $g$ traces the same boundary at least once. We conclude that:
\begin{equation}
\label{eq3}
\frac{1}{2\pi}\int_{0}^{2\pi}|f'(re^{i\theta})|d\theta\le\frac{1}{2\pi}\int_{0}^{2\pi}|g'(re^{i\theta})|d\theta.
\end{equation}
Equations (\ref{eq1}),(\ref{eq2}) and (\ref{eq3}) prove the assertion of the theorem for a general holomorphic $g$
(given that the assertion is known to be valid for one-to-one holomorphic mappings).

Thus from now on we can assume that the mapping $f$ in the statement of the theorem is one-to-one.
We solve the problem $N(2,\alpha)$ which is presented in remark 3.7. We will get the minimizing function
$\alpha z$. We clearly have the identity:
$$
\alpha^{2}=\frac{1}{2\pi}\int_{0}^{2\pi}|f(re^{i\theta})|^{2}d\theta.
$$
The perimeter of the minimal circle is $2\pi\alpha$ while the perimeter of $f(rU)$ is given by:
$$
r\int_{0}^{2\pi}|f'(re^{i\theta})|d\theta.
$$
We conclude the following inequality:
$$
2\pi\left(\frac{1}{2\pi}\int_{0}^{2\pi}|f(re^{i\theta})|d\theta\right)^{1/2}\le r\int_{0}^{2\pi}|f'(re^{i\theta})|d\theta.
$$
This proves our inequality. Finally, if we take $f(z)=\alpha z$ the inequality becomes an equality. This proves
that our inequality is, indeed sharp. $\qed $

\section{The problems $N(p,\alpha)$ for values $2<p<\infty$ of the parameter}

\begin{remark}
The two inequalities of Theorem 3.8 were proven for the value $p=2$, using results on conformal mappings. What
can be said, at this point, on similar inequalities but for values of the parameter $p$ different from $2$?
We will not tackle that problem directly. Instead, we will use convexity arguments in the form of
interpolation theory of operators. Interpolation theorems rely on two estimates given for two different
values of a parameter such as $p$, and extend them by giving estimates for all the values of that parameter
that reside between the two first values. Not always, though the inequalities for the intermediate values of
$p$ are sharp. That might happen also in cases in which the two extreme estimates are sharp. At this point
we have our inequality (which is sharp) for the value $p=2$. Lemma 3.10 below provides the second (sharp) inequality
for $p=\infty$. This case is much easier than the case $p=2$.
\end{remark}

\begin{lemma}
Let the function $f\in H(U)$ satisfy $f(0)=0$. Then for any value of $r$, $0\le r<1$ we have the following 
estimate:
$$
\max_{0\le\theta<2\pi}|f(re^{i\theta})|\le\frac{r}{2}\int_{0}^{2\pi}|f'(re^{i\theta})|d\theta.
$$
In particular we have the inequality $||f||_{\infty}\le\pi||f'||_{1}$. Both inequalities above are sharp.
\end{lemma}
\noindent
{\bf Proof.} \\
If the two inequalities of the lemma are true for $f\in H(U)$, satisfying $f(0)=0$ which are also one-to-one,
then like in the first part of the proof of Theorem  3.8 it follows that the inequalities remain true in
the more general case, where $f$ is not necessarily one-to-one. Thus we will assume from now and till the end
of the proof that the mapping $f$ in the statement of the lemma is also one-to-one. In this case we have as
in the proof of Theorem 3.8 at our disposal elementary facts from plane geometry. The expression:
$$
r\int_{0}^{2\pi}|f'(re^{i\theta})|d\theta,
$$
is the perimeter of the domain $f(rU)$. Clearly we have:
$$
2\max_{0\le\theta<2\pi}|f(re^{i\theta})|\le r\int_{0}^{2\pi}|f'(re^{i\theta})|d\theta.
$$
This last inequality is evident because $0\in f(rU)$ and because there exists a point $w\in\overline{f(rU)}$ for which:
$w=\max_{0\le\theta<2\pi}|f(re^{i\theta})|$. The first inequality of the assertion follows. By remark 1.4 and by
example 1.5 that follows it, we conclude that the inequality is, indeed sharp. $\qed $

\begin{theorem}
Let $f\in H(U)$ satisfy $f(0)=0$, and let $2\le p\le\infty$. Then for each value of $r$, $0\le r<1$ we have:
$$
\left(\frac{1}{2\pi}\int_{0}^{2\pi}|f(re^{i\theta})|^{p}d\theta\right)^{1/p}\le
\frac{r}{2\pi^{2/p}}\int_{0}^{2\pi}|f'(re^{i\theta})|d\theta.
$$
In particular also the following inequality holds true: $||f||_{p}\le\pi^{(p-2)/p}||f'||_{1}$.
\end{theorem}
\noindent
{\bf Proof.} \\
The cases $p=2,\infty$ were proved in Theorem 3.8 and in Lemma 3.10, respectively. We choose a small $\epsilon>0$
and assume that $2<p<\infty$. There exists an $N$, $p<N$ such that for any $N<q$ we have the following estimate:
\begin{equation}
\label{eq4}
\frac{(1/(2\pi))\int_{0}^{2\pi}|f(re^{i\theta})|^{q}d\theta}{((r/(2\pi))\int_{0}^{2\pi}|f'(re^{i\theta})|d\theta)^{q}}
\le\pi^{q-2+\epsilon}.
\end{equation}
This follows from Lemma 3.10. In the band $0<\Re\{z\}<1$ we define the following function:
$$
F(z)=\frac{(1/(2\pi))\int_{0}^{2\pi}|f(re^{i\theta})|^{2z+q(1-z)}d\theta}
{((r/(2\pi))\int_{0}^{2\pi}|f'(re^{i\theta})|d\theta)^{2z+q(1-z)}}.
$$
Then $F(z)$ is analytic in $0<\Re\{z\}\le 1$ (for if $c>0$, then $c^{z}$ is non-zero analytic).
Let us write $z=t+is$, $t,s\in\mathbb{R}$, $0\le t\le 1$. By the triangle inequality we obtain:
$$
|F(z)|\le\frac{(1/(2\pi))\int_{0}^{2\pi}|f(re^{i\theta})|^{2t+q(1-t)}d\theta}
{((r/(2\pi))\int_{0}^{2\pi}|f'(re^{i\theta})|d\theta)^{2t+q(1-t)}}.
$$
Hence: \\
by equation (\ref{eq4}) we have $|F(is)|\le\pi^{q-2+\epsilon}$, and \\
by Theorem 3.8 $|F(1+is)|\le 1$. \\
By Hadamard Convexity Theorem we deduce that:
\begin{equation}
\label{eq5}
|F(t+is)|\le\pi^{(q-2+\epsilon)(1-t)}.
\end{equation}
Now, suppose that $p=2t+q(1-t)$, then:
$$
(q-2+\epsilon)(1-t)=(p-2)+\epsilon\left(\frac{p-2}{q-2}\right).
$$
So by equation (\ref{eq5}) we get:
$$
\left(\frac{1}{2\pi}\int_{0}^{2\pi}|f(re^{i\theta})|^{p}d\theta\right)^{1/p}\le
\left(\frac{r}{2\pi^{2/p}}\int_{0}^{2\pi}|f'(re^{i\theta})|d\theta\right)\pi^{\epsilon(p-2)/(q-2)}.
$$
Now, we will take $q\rightarrow\infty$ and obtain:
$$
\left(\frac{1}{2\pi}\int_{0}^{2\pi}|f(re^{i\theta})|^{p}d\theta\right)^{1/p}\le
\frac{r}{2\pi^{2/p}}\int_{0}^{2\pi}|f'(re^{i\theta})|d\theta.
$$
The proof of the theorem is now complete. $\qed $ \\

\begin{theorem}
If $f\in H(U)$ and if $||f'||_{1}<\infty$, then $f\in H^{p}(U)$ for all $p$, $2\le p\le\infty$ and
the following inequality holds true:
$$
||f||_{p}\le\pi^{(p-2)/p}||f'||_{1}+|f(0)|.
$$
\end{theorem}
\noindent
{\bf Proof.} \\
We use Theorem 4.3 with the function $f(z)-f(0)$, and than use the triangle inequality. $\qed $ \\

\begin{remark}
The inequalities in Theorem 4.3 are sharp for $p=2,\infty$ and we get equalities for the extremal
function $f(z)=\alpha z$, in both cases. We do not know if these inequalities are sharp for the
other values of $p$, i.e. $2<p<\infty$. In any event for these values of the parameter $p$, the 
function $f(z)=\alpha z$ does not give us equality. Thus either the inequalities in Theorem 4.3 
are not sharp, or they are sharp but for values $2<p<\infty$, the function $f(z)=\alpha z$ is not 
an extremal function. 
\end{remark}

\section{Few open problems}

Here are four natural problems for which we do not know the answer at this point: \\
\\
1. Let $f\in S$ and let $F$ be the Steiner symmetrization of $f$. For which values of $p$,
$2\le p\le\infty$ we have the following inequality?
$$
\int_{0}^{2\pi}|f(re^{i\theta})|^{p}d\theta\le\int_{0}^{2\pi}|F(re^{i\theta})|^{p}d\theta.
$$
We know that this inequality holds true for $p=2$ (by Theorem 2.7) and it is faulty for $p=\infty$ (easy). \\
\\
2. Let $f\in S(2,\alpha)$ and let $\{c_{\phi_n}\}_{n=1}^{\infty}$ be a sequence of shrinking factors
that deform $f(U)$ into the minimal disk. What is the value of the convergent infinite product
$\prod_{n=1}^{\infty}c_{\phi_n}$? \\
\\
3. Are the inequalities in Theorem 4.3 and in Theorem 4.4 sharp? If not find the optimal constants
(which should be at most $\pi^{(p-2)/p}$) and find the optimal functions. \\
\\
4. Are the extremal functions for the inequalities in Theorem 4.3 (with optimal constants) unique?

\section{Inequalities on the real part of Steiner symmetrization}

The results of this section might be already known to experts on Steiner symmetrization. The author
was not able to find in the literature these results. We refer to the first sentence in section 4.11
on page 130 of the book \cite{ha}. Here it is: "Circular symmetrization is more powerful than Steiner
symmetrization, and any result obtainable by the later method can also be obtained by the former on
taking exponentials, though this may be less direct." The author did not check the validity of this
declaration, however, the hint of connecting the two types of symmetrization by the exponential function
was taken in order to try and resolve problem number 1 on the list of problems given in the previous section (section 5).
We did not manage to solve that problem. However, this idea in W. K. Hayman's book,\cite{ha}, produced
a family of integral inequalities comparing $\Re\{f(z)\}$ with the real part of the Steiner symmetrization
$\Re\{F(z)\}$ of $f(z)$. We remark that connecting the Circular symmetrization with the Steiner symmetrization
by the exponential function is a very different idea than the geometric idea of connecting them by
shifting the function to infinity. This idea was presented in the proof of Theorem 2.7 and gave us
a positive answer to problem 1 above, in the case where $p=2$.

Let $f\in H(U)$ and let us denote $D=f(U)$. Let $D^{*}$ be the Steiner symmetrization of $D$, with respect to
the $x$-axis. Then $D^{*}$ is a simply connected domain. We recall that for each $a\in\mathbb{R}$, we denote:
$l(a)={\rm meas}\{y\in\mathbb{R}\,|\,a+iy\in D\}$, and we have by the definition of the Steiner symmetrization:
$D^{*}=\{a+iy\,|\,|y|< (1/2)\cdot l(a)\}$. This means that for any $a\in\mathbb{R}$ the intersection of the
vertical line $x=a$ with $D^{*}$ is either an empty set or an open vertical interval (line segment) symmetric
about the $x$-axis at the point of the intersection, $(a,0)$: $\{x=a\}\cap D^{*}=\{a+iy\,|\,|y|<(1/2)\cdot l(a)\}$.
This vertical interval might be a full line parallel to the $y$-axis, in the case that $l(a)=\infty$. 

Let us apply (as suggested by Hayman) the exponential mapping to this line segment:
\begin{equation}
\label{eq6}
\exp\left(\{x=a\}\cap D^{*}\right)=\left\{e^{a}e^{iy}\,\left|\,|y|<\frac{1}{2}l(a)\right.\right\}.
\end{equation}

This is a circular arc. This arc is a proper arc of the circle $|z|=e^{a}$ provided that $l(a)\le 2\pi$.
If $2\pi<l(a)$ then the arc visits some of the points of the circle $|z|=e^{a}$ at least twice. This means that
in any event the set $\exp(D^{*})$ is a circular symmetric domain. If for all $a\in\mathbb{R}$ we
have $l(a)<2\pi$ then the set $\exp(D^{*})$ is a simply connected circular symmetric domain. For any
$a\in\mathbb{R}$ we have the following identity:
$$
\{|z|=e^{a}\}\cap\exp(D)=\{z\,|\,|z|=e^{w}, w\in D, \Re\{w\}=a\}.
$$
We note that the total length of the arcs that comprise the set $\{|z|=e^{a}\}\cap\exp(D)$ is $e^{a}\cdot l(a)$,
so when we form the P\"olya symmetrization of $\exp(D)$ we obtain:
$$
\bigcup_{a\in\mathbb{R}}\left\{a^{a}\cdot e^{i\theta}\,\left|\,|\theta|<\frac{1}{2}l(a)\right.\right\},
$$
and by equation (\ref{eq6}) this is:
$$
\bigcup_{a\in\mathbb{R}}\exp(\{x=a\}\cap D^{*})=\exp(D^{*}).
$$
This proves parts 1 and 2 of the following:

\begin{proposition}
1. $\exp(D^{*})$ is the circular symmetrization of $\exp(D)$, where $D$ is a domain and $D^{*}$ is the
Steiner symmetrization of this domain. \\
2. If for all $a\in\mathbb{R}$, $l(a)<2\pi$, then the intersection arcs $\{e^{a}e^{i\theta}\,|\,|\theta|<(1/2)l(a)\}$
are simple (i.e. they do not pass through any point more than once). \\
3. If for all $a\in\mathbb{R}$, $l(a)<2\pi$, then $\exp(D^{*})$ is a simply connected domain.
\end{proposition}
\noindent
{\bf Proof.} \\
We only need to prove part 3 (because part 1 and part 2 were proven above). However, part 3 follows at once by 
part 2 and by the fact that the Steiner symmetric domain $D^{*}$ is a simply connected domain. $\qed $ \\
\\
Using the result of Albert Baernstein (Theorem 6 in \cite{b}) we have the following:

\begin{theorem}
Let $f\in H(U)$ and let us denote $D=f(U)$, and assume that for any $a\in\mathbb{R}$ we have $l(a)<2\pi$.
Also suppose that $f(0)\ge 0$. Let $F\in H(U)$ be a conformal mapping of $U$ onto $D^{*}$, where
$F(0)=|f(0)|=f(0)$ and where $D^{*}$ is the Steiner symmetrization of $D$. \\
If $\Phi$ is a convex non-decreasing function on $(-\infty,\infty)$, then for all $r$, $0\le r<1$, we have:
$$
\int_{-\pi}^{\pi} \Phi(\Re\{f(re^{i\theta})\})d\theta\le\int_{-\pi}^{\pi} \Phi(\Re\{F(re^{i\theta})\})d\theta.
$$
\end{theorem}
\noindent
{\bf Proof.} \\
By Proposition 6.1(1): $\exp(D^{*})$ is the circular symmetrization of $\exp(D)$. By our assumption on
$f$ ($\forall\,a\in\mathbb{R}$, $l(a)<2\pi$), and by Proposition 6.1(2), the mapping: $\exp\,:\,D^{*}\rightarrow
\exp(D^{*})$, is injective and so conformal. We recall that the mapping: $F\,:\,U\rightarrow D^{*}$ is
conformal and $F(0)=|f(0)|$. Hence the composition: $\exp(F)\,:\,U\rightarrow\exp(D^{*})$, is a conformal
mapping of $U$ onto $\exp(D^{*})$ which is a simply connected domain (as should be the case), by Proposition 6.1(3).
Also we have: $\exp(F(0))=\exp(|f(0)|)=\exp(f(0))$, where the last equality follows by our assumption, $f(0)\ge 0$.

To sum up we have $g(z)=\exp(f(z))\in H(U)$ where by the above notations: $g\,:\,U\rightarrow \exp(D)=g(U)$. The
mapping $G(z)=\exp(F(z))$ is a conformal and onto mapping $G\,:\,U\rightarrow\exp(D^{*})$ that satisfies
$G(0)=\exp(F(0))=\exp(f(0))=g(0)=|g(0)|$. The simply connected domain $\exp(D^{*})$ is the circular
symmetrization of $\exp(D)=g(U)$. Thus the pair of mappings $g$, $G$ satisfy all the assumptions of Albert
Baernstein result, Theorem 6 in \cite{b}. Using this theorem we obtain:
$$
\int_{-\pi}^{\pi}\Phi(\log|g(re^{i\theta})|)d\theta\le\int_{-\pi}^{\pi}\Phi(\log|G(re^{i\theta})|)d\theta,
$$
for any convex and non-decreasing $\Phi$ on $(-\infty,\infty)$, and any $r$, $0\le r<1$. Plugging in
the expressions $g(z)=\exp(f(z))$ and $G(z)=\exp(F(z))$ we obtain:
$$
\int_{-\pi}^{\pi}\Phi(\log|\exp(f(re^{i\theta}))|)d\theta\le\int_{-\pi}^{\pi}\Phi(\log|\exp(F(re^{i\theta}))|)d\theta.
$$
Finally, since $\log|\exp(\alpha)|=\Re\{\alpha\}$ for any $\alpha\in\mathbb{C}$, we get:
$$
\int_{-\pi}^{\pi} \Phi(\Re\{f(re^{i\theta})\})d\theta\le \int_{-\pi}^{\pi} \Phi(\Re\{F(re^{i\theta})\})d\theta,
$$ 
where $F$ is the Steiner symmetrization of $f$. $\qed $ \\
\\
If we take $\Phi(x)=e^{pX}$ for some $p>0$, and take,
$$
\Phi(x)=\left\{\begin{array}{lll} x^{p} & , & x\ge 0 \\ 0 & , & x<0 \end{array}\right.
$$
for some $p>1$, then we note that both functions are convex non-decreasing on $(-\infty,\infty)$, and we deduce
from Theorem 6.2 the following:

\begin{corollary}
Let $f\in H(U)$ satisfy $l(a)<2\pi$ for all $a\in\mathbb{R}$ and $f(0)\ge 0$, and let $F\,:\,U\rightarrow f(U)^{*}$
($f(U)^{*}$ is the Steiner symmetrization of $f(U)$) be a conformal onto with $F(0)=f(0)$, then: \\
1. For any $0<p$ we have:
$$
\int_{-\pi}^{\pi}\exp\left(p\Re\{f(re^{i\theta})\}\right)d\theta\le\int_{-\pi}^{\pi}\exp\left(p\Re\{F(re^{i\theta})\}\right)d\theta,\,\,\,
0\le r<1.
$$
2. For any $1<p$ we have:
$$
\int_{-\pi}^{\pi}\left(\Re\{f(re^{i\theta})\}\right)_{+}^{p}d\theta\le\int_{-\pi}^{\pi}\left(\Re\{F(re^{i\theta})\}\right)_{+}^{p}d\theta,
\,\,\,0\le r<1.
$$
Here if $a\in\mathbb{R}$, then we denote:
$$
a_{+}=\left\{\begin{array}{lll} a & , & 0\le a \\ 0 & , & a<0\end{array}\right.
$$
In particular if in 2 we take the $p$'th root from both sides and than let $p\rightarrow\infty$, we obtain:
$$
\sup_{\theta}\left(\Re\{f(re^{i\theta})\}\right)_{+}\le\sup_{\theta}\left(\Re\{F(re^{i\theta})\}\right)_{+},\,\,\,0\le r<1.
$$
\end{corollary}
The last inequality in Corollary 6.3 is clear by the definition of the Steiner symmetrization, for the value
$r=1$ (if that makes sense).

\begin{remark}
1. In fact it follows by the definition of Steiner symmetrization that: if $h(z)=\Re\{f(z)\}$ and if $H(z)=\Re\{F(z)\}$,
then we have $h(U)=H(U)$ and in particular:
$$
\lim_{r\rightarrow 1^{-}}\inf_{\theta}\Re\{f(re^{i\theta})\}=\lim_{r\rightarrow 1^{-}}\inf_{\theta}\Re\{F(re^{i\theta})\},
$$
and also
$$
\lim_{r\rightarrow 1^{-}}\sup_{\theta}\Re\{f(re^{i\theta})\}=\lim_{r\rightarrow 1^{-}}\sup_{\theta}\Re\{F(re^{i\theta})\}.
$$
2. Another example for a concrete inequality we can deduce from Theorem 6.2 is, for example, the following: we take for
$\Phi(x)$ the function $\Phi(x)=\exp(x_{+}^{p})$, $p>1$. The first two derivatives on $x>0$ are,
$$
\Phi'(x)=px^{p-1}\exp(x^{p})>0,\,\,\,\Phi''(x)=px^{p-2}(p-1+px^{p})\exp(x^{p})>0.
$$
For $x\le 0$ we have $\Phi(x)\equiv 1$, a constant function. Using this $\Phi(x)$, we obtain by Theorem 6.2 the following 
result:
\end{remark}

\begin{corollary}
Let $f\in H(U)$ satisfy $l(a)<2\pi$ for all $a\in\mathbb{R}$, and $f(0)\ge 0$, and let $F\,:\,U\rightarrow f(U)^{*}$ be
conformal, onto with $F(0)=f(0)$, then for any $p>1$ we have the following inequality:
$$
\int_{-\pi}^{\pi}\exp\left(\left(\Re\{f(re^{i\theta})\}\right)_{+}^{p}\right)d\theta\le
\int_{-\pi}^{\pi}\exp\left(\left(\Re\{F(re^{i\theta})\}\right)_{+}^{p}\right)d\theta,\,\,\,{\rm for}\,{\rm all}\,\,0\le r<1.
$$
\end{corollary}

\begin{remark}
In fact Corollary 6.5 is a special case of the following more general statement: \\
If $f(z)$ and $F(z)$ are as in Corollary 6.5, and if $\Psi(z)$ is any entire function with non-negative MaClaurin
coefficients (i.e. $\Psi^{(n)}(0)\ge 0$ for any $n\in\mathbb{Z}^{+}\cup\{0\}$, then for any $p>1$ and any $0\le r<1$
we have:
$$
\int_{-\pi}^{\pi}\Psi\left(\left(\Re\{f(re^{i\theta})\}\right)_{+}^{p}\right)d\theta\le
\int_{-\pi}^{\pi}\Psi\left(\left(\Re\{F(re^{i\theta})\}\right)_{+}^{p}\right)d\theta.
$$
\noindent
{\bf Proof.} \\
By the assumptions we have: $\Psi(z)=\sum_{n=0}^{\infty}a_{n}z^{n}$, with $a_{n}\ge 0$ for $n\in\mathbb{Z}^{+}\cup\{0\}$.
Now the claim follows by Corollary 6.3(2). For by that corollary and the non-negativity of the coefficients $a_{n}$, we get:
$$
\int_{-\pi}^{\pi}a_{n}\left(\Re\{f(re^{i\theta})\}\right)_{+}^{pn}d\theta\le
\int_{-\pi}^{\pi}a_{n}\left(\Re\{F(re^{i\theta})\}\right)_{+}^{pn}d\theta,
$$
$$
\int_{-\pi}^{\pi}\sum_{n=0}^{\infty}a_{n}\left(\left(\Re\{f(re^{i\theta})\}\right)_{+}^{p}\right)^{n}d\theta\le
\int_{-\pi}^{\pi}\sum_{n=0}^{\infty}a_{n}\left(\left(\Re\{F(re^{i\theta})\}\right)_{+}^{p}\right)^{n}d\theta,
$$
$$
\int_{-\pi}^{\pi}\Psi\left(\left(\Re\{f(re^{i\theta})\}\right)_{+}^{p}\right)d\theta\le
\int_{-\pi}^{\pi}\Psi\left(\left(\Re\{F(re^{i\theta})\}\right)_{+}^{p}\right)d\theta.
$$
$\qed $ \\
\\
The change of order in the integration and the summation is easily justified.
\end{remark}

\begin{remark}
In the above results, the function $f\in H(U)$ was assumed to satisfy $l(a)<2\pi$, $\forall\,a\in\mathbb{R}$. In the case
we had another uniform (finite) upper bound, say $l(a)<M$, $\forall\,a\in\mathbb{R}$, we could have looked (in case $2\pi<M$)
in the scaled function, $2\pi f(z)/M$, and get instead of the inequality of Theorem 6.2 the following inequality:
$$
\int_{-\pi}^{\pi} \Phi(\frac{2\pi}{M}\Re\{f(re^{i\theta})\})d\theta\le
\int_{-\pi}^{\pi} \Phi(\frac{2\pi}{M}\Re\{F(re^{i\theta})\})d\theta.
$$
What if there was no uniform upper bound on the $l(a)$'s, i.e. $\sup_{a\in\mathbb{R}} l(a)=+\infty$? We would not like to conformaly
map $f(U)$ into a bounded subset of $\mathbb{C}$, say by an inversion:
$$
\left(\frac{\alpha}{z-\beta}\right).
$$
The reason is that there does not seem to be a simple relation between the Steiner symmetrization of the original function
$f(z)$ and the Steiner symmetrization of the transformed mapping:
$$
\left(\frac{\alpha}{f(z)-\beta}\right).
$$
What we can do is the following: Pick a number $r_{0}$, $0<r_{0}<1$ and consider the function $f(r_{0}z)$. We denote by
$F(r_{0},z)$ the Steiner symmetrization of $f(r_{0}z)$. Since the image: $\{f(r_{0}z)\,|\,|z|<1\}$ is a bounded set, our
theorems give us comparison between integrals that involve $f(r_{0}z)$ and those that involve $F(r_{0},z)$. We now have to
estimate the relations between integrals that involve $F(r_{0},z)$ and those that involve the Steiner symmetrization $F(re^{i\theta})$,
when $r_{0}\rightarrow 1^{-}$. For example, is the following limit claim holds true? $\lim_{r_{0}\rightarrow 1^{-}}F(r_{0},z)=F(z)$,
$|z|<1$. In what sense? (uniform, uniform on compacta, other). This might not be easy, for we treating the case in which
$f(U)$ is unbounded vertically, i.e. $\sup_{a\in\mathbb{R}} l(a)=+\infty$ while $f(r_{0}U)$ is a bounded set for $0<r_{0}<1$.
The mappings $F(z)\,:\,U\rightarrow f(U)^{*}$, $F(r_{0},z)\,:\,U\rightarrow f(r_{0}U)^{*}$ are conformal and onto and satisfy
the conditions: $F(0)=f(0)=f(r_{0}\cdot 0)=F(r_{0},0)$.

We need a kind of a continuity claim on families of conformal mappings $U\rightarrow D_{r_{0}}$ such that $D_{r_{0}}\rightarrow D$
when $r_{0}\rightarrow 1^{-}$ in some sense (What do we mean by $D_{r_{0}}\rightarrow D$?). Once again, in our model the
domains $D_{r_{0}}$ ($0<r_{0}<1$) are bounded, while $D$ is unbounded.
\end{remark}

\section{Steiner symmetrization and zero sets of bounded holomorphic functions in $U$}

In this section we will point at some facts related to Problem number 2, on the list of problems in section 5. We will use
the following definition: for $2\le p\le\infty$, $0<\alpha<\infty$ we define $S(p,\alpha)=\{f\in H(U)\,|\,f\,\,{\rm is}\,{\rm univalent}\,
{\rm in}\, U, f(0)=0, 1\le f'(0), \alpha\le ||f||_{p}\}$. We note that our normalization in the defining equation of 
$S(p,\alpha)$ included the inequality, $1\le |f'(0)|$ while now we gave up the absolute value and use instead $1\le f'(0)$.
Problem 2 asks the following: Let $f\in S(2,\alpha)$ and let
$\{c_{\phi_{n}}\}_{n=1}^{\infty}$ be a sequence of shrinking factors that deform $f(U)$ into the minimal disk. What is
the value of the following convergent infinite product $\prod_{n=1}^{\infty}c_{\phi_{n}}$?

For $f\in S(2,\alpha)$ and $\phi\in\mathbb{R}$, we denoted by $g_{\phi}$ the Steiner symmetrization of $e^{i\phi}f(z)$.
It has the following four properties: \\
\\
1. $g_{\phi}\in H(U)$ is one-to-one in $U$, and $g_{\phi}(0)=0$. \\
\\
2. $1\le |e^{i\phi}f'(0)|\le g_{\phi}^{'}(0)$. Note that $|e^{i\phi}f'(0)|=f'(0)\ge 1$. \\
\\
3. $\alpha\le ||e^{i\phi}f||_{2}\le ||g_{\phi}||_{2}$. Note that $||e^{i\phi}f||_{2}=||f||_{2}\ge\alpha$. \\
\\
4. $\int_{-\pi}^{\pi}|g_{\phi}^{'}(e^{i\theta})|d\theta\le\int_{-\pi}^{\pi}|e^{i\phi}f'(e^{i\theta})|d\theta$. Note that
$\int_{-\pi}^{\pi}|e^{i\phi}f'(e^{i\theta})|d\theta=\int_{-\pi}^{\pi}|f'(e^{i\theta})|d\theta$. \\
\\
The shrinking factor that corresponds to $f$ and to $\phi$ is the following number:
$$
c_{\phi}=\max\left\{\frac{1}{g_{\phi}^{'}(0)},\frac{\alpha}{||g_{\phi}||_{2}}\right\}.
$$
Thus by properties 2 and 3 we have $0< c_{\phi}\le 1$. By Theorem 3.4, if $0<\alpha<\infty$, and if $f\in S(2,\alpha)$ and
$\{\phi_{n}\}_{n=1}^{\infty}$ is any sequence of real numbers, then $0<\prod_{n=1}^{\infty}c_{\phi_{n}}\le 1$. This theorem
was the source of Problem 2.

\begin{proposition}
\begin{equation}
\label{eq7}
\sum_{n=1}^{\infty}\log\left(\frac{1}{c_{\phi_{n}}}\right)<\infty.
\end{equation}
\begin{equation}
\label{eq8}
\sum_{n=1}^{\infty}(1-c_{\phi_{n}})<\infty.
\end{equation}
\end{proposition}
\noindent
{\bf Proof.} \\
We prove equation (\ref{eq7}): By Theorem 3.4 we have: $\log\prod_{n=1}^{\infty}c_{\phi_{n}}=
\sum_{n=1}^{\infty} \log c_{\phi_{n}}>-\infty$. Hence $\sum_{n=1}^{\infty}\log (1/c_{\phi_{n}})<\infty$. \\
Now we prove equation (\ref{eq8}): We write $c_{\phi_{n}}=1-b_{\phi_{n}}$. Then $0\le b_{\phi_{n}}<1$. Also we
have $\prod_{n=1}^{\infty}c_{\phi_{n}}=\prod_{n=1}^{\infty}(1-b_{\phi_{n}})$. For $0\le b<1$ we have 
$e^{-b}=1-b+b^{2}/2!-\ldots\ge 1-b$, because this MaClaurin expansion is a Leibniz series. This also implies that:
$0\le e^{-b}-(1-b)\le b^{2}/2$. So $\prod_{n=1}^{\infty}e^{b_{\phi_{n}}}\ge\prod_{n=1}^{\infty}(1-b_{\phi_{n}})>0$.
Thus: $\exp(-\sum_{n=1}^{\infty}b_{\phi_{n}})>0$ and so $\sum_{n=1}^{\infty} b_{\phi_{n}}<\infty$. But $b_{\phi_{n}}=
1-c_{\phi_{n}}$ and we conclude that $\sum_{n=1}^{\infty}(1-c_{\phi_{n}})<\infty$. $\qed $ \\
\\
The inequality $\sum_{n=1}^{\infty} (1-c_{\phi_{n}})<\infty$ says that the sequence $\{c_{\phi_{n}}\}_{n=1}^{\infty}$
satisfies the Blaschke condition. Hence this sequence is precisely the zero set of a bounded analytic function in
$U$. In fact the corresponding Blaschke product converges in $U$:
$$
\prod_{n=1}^{\infty}\left(\frac{z-c_{\phi_{n}}}{1-c_{\phi_{n}}z}\right).
$$

\begin{proposition}
If $0<\alpha<\infty$, and if $f\in S(2,\alpha)$, and if $\{\phi_{n}\}_{n=1}^{\infty}$ is any sequence
of real numbers, then the infinite product:
$$
B_{\{\phi_{n}\}}(z)=\prod_{n=1}^{\infty}\left(\frac{z-c_{\phi_{n}}}{1-c_{\phi_{n}}z}\right),
$$
is a Blaschke product, i.e. it is uniformly convergent on compact subsets of $U$, and $\{c_{\phi_{n}}\}_{n=1}^{\infty}$ is
the zero set of the resulting bounded (by $1$) analytic function, $B_{\{\phi_{n}\}}(z)$.
\end{proposition}

This naturally leads to the question: is the converse of Proposition 7.2 holds true? \\
\\
\underline{Problem 2'.} Let:
$$
\prod_{n=1}^{\infty}\left(\frac{z-c_{\alpha_{n}}}{1-c_{\alpha_{n}}z}\right)
$$
be a Blaschke product all of whose zeros $\{\alpha_{n}\}_{n=1}^{\infty}$ are positive numbers, $0<\alpha_{n}\le 1$.
Is there a number $0<\alpha<\infty$ and a function $f\in S(2,\alpha)$, and a sequence of real numbers 
$\{\phi_{n}\}_{n=1}^{\infty}$ such that $\forall\,n\in\mathbb{Z}^{+}$, $c_{\phi_{n}}=\alpha_{n}$ the 
corresponding shrinking factors? \\
\\
\underline{An explanation.} We are given the data $\{\alpha_{n}\}_{n=1}^{\infty}$ and should come up
with an $0<\alpha<\infty$, $f\in S(2,\alpha)$ and real numbers $\{\phi_{n}\}_{n=1}^{\infty}$ such that: \\
1. If $\overline{g}_{\phi_{1}}$ is the Steiner symmetrization of $e^{i\phi_{1}}f$ then,
$$
c_{\phi_{1}}=\max\left\{\frac{1}{\overline{g}_{\phi_{1}}'(0)},\frac{\alpha}{||\overline{g}_{\phi_{1}}||_{2}}\right\}
=\alpha_{1},\,\,\,{\rm and}\,\,g_{\phi_{1}}=c_{\phi_{1}}\overline{g}_{\phi_{1}}.
$$
2. If $(\overline{g_{\phi_{1}}})_{\phi_{2}}$ is the Steiner symmetrization of $e^{i\phi_{2}}g_{\phi_{1}}$ then,
$$
c_{\phi_{2}}=\max\left\{\frac{1}{\overline{g_{\phi_{1}}}_{\phi_{2}}'(0)},
\frac{\alpha}{||(\overline{g_{\phi_{1}}})_{\phi_{2}}||_{2}}\right\}
=\alpha_{2},\,\,\,{\rm and}\,\,(g_{\phi_{1}})_{\phi_{2}}=c_{\phi_{2}}\overline{(g_{\phi_{1}}})_{\phi_{2}}.
$$
3. If $(\overline{(g_{\phi_{1}})_{\phi_{2}}})_{\phi_{3}}$ is the Steiner symmetrization of 
$e^{i\phi_{3}}(g_{\phi_{1}})_{\phi_{2}}$ then,
$$
c_{\phi_{3}}=\max\left\{\frac{1}{(\overline{(g_{\phi_{1}})_{\phi_{2}}})_{\phi_{3}}'(0)},
\frac{\alpha}{||(\overline{(g_{\phi_{1}})_{\phi_{2}}})_{\phi_{3}}||_{2}}\right\}
=\alpha_{3},\,\,\,{\rm and}\,\,((g_{\phi_{1}})_{\phi_{2}})_{\phi_{3}}=
c_{\phi_{3}}(\overline{(g_{\phi_{1}})_{\phi_{2}}})_{\phi_{3}},\,\,{\rm etc}\ldots
$$
\underline{An idea.} \\
We show that for each $N\in\mathbb{Z}^{+}$ we can construct a number $0<\alpha<\infty$, a function $f_{N}\in S(2,\alpha)$
and real numbers $\{\phi_{n}^{(N)}\}_{n=1}^{\infty}$ such that this solves the finite problem, i.e. $\forall\,n\in\mathbb{Z}^{+}$,
$1\le n\le N$, $c_{\phi_{n}^{(N)}}=\alpha_{n}$ the corresponding shrinking factors. This gives a sequence of functions
$\{f_{N}\}_{N=1}^{\infty}\subseteq S(2,\alpha)$, and a sequence of sequences of shrinking factors $\{c_{\phi_{n}^{(N)}}\}_{n=1}^{N}=
\{\alpha_{1},\ldots,\alpha_{N}\}$ on the numbers $\{\phi_{n}^{(N)}\}_{n=1}^{\infty}$. We then might try to prove that
the Blaschke condition $\sum_{n=1}^{\infty}(1-\alpha_{n})<\infty$ implies that the limit $f(z)=\lim{N\rightarrow\infty}f_{N}(z)$
exists uniformly on compact subsets of $U$. We also might hope to prove that the infinite set of limits:
$\lim_{N\rightarrow\infty}\phi_{n}^{(N)}=\phi_{n}^{(\infty)}$ exist and finally that the data $f(z)\in S(2,\alpha)$ and
$\{\phi_{n}^{(\infty)}\}_{n=1}^{\infty}$ solves Problem 2', namely that we have $c_{\phi_{n}^{(\infty)}}=\alpha_{n}$,
$\forall\,n\in\mathbb{Z}^{+}$. So we first want to solve: \\
\\
\underline{Problem 2 (finite).} Given $N$ numbers $\alpha_{n}$, $0<\alpha_{1}<\alpha_{2}<\ldots<\alpha_{N}<1$, find a
number $0<\alpha<\infty$ and a function $f_{N}\in S(2,\alpha)$ and a sequence of real numbers $\{\phi_{n}\}_{n=1}^{N}$
such that $c_{\phi_{n}}=\alpha_{n}$, $n=1,2,\ldots,N$.

\begin{remark}
We recall that we have defined the shrinking factor as follows:
$$
c_{\phi}=\max\left\{\frac{1}{g_{\phi}^{'}(0)},\frac{\alpha}{||g_{\phi}||_{2}}\right\},\,\,{\rm so}\,\,0<c_{\phi}\le 1,
$$
and since $c_{\phi}=1-b_{\phi}$, also $b_{\phi}=1-c_{\phi}$. Thus we have the following identities:
$$
b_{\phi}=\min\left\{1-\frac{1}{g_{\phi}^{'}(0)},1-\frac{\alpha}{||g_{\phi}||_{2}}\right\},
$$
and
$$
\frac{1}{c_{\phi}}=\min\left\{g_{\phi}^{'}(0),\frac{1}{\alpha}\cdot ||g_{\phi}||_{2}\right\}.
$$
Thus we can restate Proposition 7.1 as follows:
\end{remark}

\begin{proposition}
\begin{equation}
\label{eq9}
\sum_{n=1}^{\infty}\log\left(\min\left\{g_{\phi_{n}}^{'}(0),\frac{1}{\alpha}\cdot||g_{\phi_{n}}||_{2}\right\}\right)<\infty.
\end{equation}
\begin{equation}
\label{eq10}
\sum_{n=1}^{\infty}\min\left\{1-\frac{1}{g_{\phi_{n}}^{'}(0)},1-\frac{\alpha}{||g_{\phi_{n}}||_{2}}\right\}<\infty.
\end{equation}
\end{proposition}
We will need to use results on the convergence of sequences of conformal mappings. The next section surveys the
results we will be using.

\section{Convergence of a sequence of conformal mappings}
We will use the following two references: \\
1. \cite{c}, sections 120 through 124, pages 74-77. \\
2. \cite{g}, section 5 pages 54-62. \\
\\
The exposition in \cite{g} is easier for us being more modern but the results were proven by Constantine Carath\'eodory
in \cite{c}. So we take the parts we need mostly from \cite{g}. Section 5 in \cite{g} is titled: "Convergence theorems
on the conformal mapping of a sequence of domains": \\
Suppose we have a sequence of univalent domains $B_{1},B_{2},\ldots$, in the $z$-plane, each including $z=0$. If there
exists a disk $|z|<\rho$, where $\rho>0$, that belongs to all the domains in $B_{n}$, we define \underline{the kernel} of
this sequence of domains as the largest domain containing $z=0$ such that an arbitrary closed subset of it belongs to all
the domains $B_{n}$ from some $n$ on. By "largest domain" is meant the domain containing any other domain possessing this
property. If such a disk does not exist, the kernel of the sequence of domains $B_{1},B_{2},\ldots$ is defined to
be the point $z=0$. We shall say that the sequence of domains $B_{1},B_{2},\ldots$ \underline{converges to the kernel} $B$,
and we shall denote this by writing $B_{n}\rightarrow B$, if every subsequence of these domains has $B$ as its kernel. In particular,
if a sequence of simply connected domains, $B_{1},B_{2},\ldots,B_{n},\ldots$ that include $z=0$ converges to the
limiting domain $B$ (also including $z=0$) in the sense that all boundary points of the domains $B_{n}$ from some $n$ 
on are arbitrary close to the boundary of the domain $B$, and all points of the boundary of the domain $B$ are arbitrary
close to the boundaries of the domains $B_{n}$, then this sequence has the domain $B$ as its kernel and it converges
to that kernel.

{\bf In our application later on the domains in the sequence of domains are Steiner symmetric and will
turn out to satisfy exactly the assumptions of the previous paragraph}.

Convergence to the kernel is guaranteed also for a sequence of domains $B_{n}$ that include $z=0$ and satisfy the
condition $B_{1}\subseteq B_{2}\subseteq B_{3}\subseteq\ldots$, or for a sequence of domains $B_{n}$ that contains a
neighborhood of the point $z=0$ and satisfy the condition $B_{1}\supseteq B_{2}\supseteq B_{3}\supseteq\ldots$. \\
\\
{\bf Theorem A (Carath\'eodory, \cite{c12}).} Suppose that we have a sequence of functions $z=f_{n}(\xi)$, where
$n=1,2,\ldots$, that are regular in the disk $|\xi|<1$. Suppose $f_{n}(0)=0$ and $f_{n}^{'}(0)>0$ for $n=1,2,\ldots$.
Suppose that, for each $n$, the function $f_{n}(\xi)$ maps the disk $|\xi|<1$ onto a domain $B_{n}$. For the sequence
$\{f_{n}(\xi)\}$ to converge in $|\xi|<1$ to a finite function, it is necessary and sufficient that the sequence 
$\{B_{n}\}$ converge to the kernel $B$, which is either the point $z=0$ or a domain having more than one boundary
point. When convergence exists, it is uniform inside the disk $|\xi|<1$. If the limit function $f(\xi)\not\equiv {\rm const.}$,
it maps $|\xi|<1$ onto the kernel $B$, and the sequence $\{\phi_{n}(z)\}$ \underline{of inverse functions} $\phi_{n}(z)$
converges uniformly inside $B$ to the function $\phi(z)$ inverse to $f(\xi)$. (Thus it is assumed that the functions 
$f_{n}(\xi)$ are \underline{conformal}). \\

\begin{remark}
It was proved in \cite{bi} and in \cite{m} that the conditions of Theorem A are also necessary and sufficient for 
convergence in mean of $\{_{n}^{'}(\xi)\}$ to $f^{'}(\xi)$, that is, necessary and sufficient for:
$$
\lim_{n\rightarrow\infty}\iint_{B}|f_{n}^{'}(\xi)-f^{'}(\xi)|^{2}d\sigma=0,
$$
where $f_{n}^{'}(\xi)$ is taken equal to $0$ outside the domain $B_{n}$
\end{remark}
Theorem A gives the conditions for convergence of univalent functions only in the open disk $|\xi|<1$. For the
convergence of univalent functions in the closed disk $|\xi|\le 1$, we give the following theorem, confining
ourselves to domains of the Jordan type. \\
\\
{\bf Theorem B (Rad\'o, \cite{r}).} Let $\{B_{n}\}$, $n=1,2,\ldots$ denote a sequence of simply connected domains
each including the point $z=0$ and each bounded by a Jordan curve. Denote the boundary of $B_{n}$ by $C_{n}$. Suppose
that the sequence $\{B_{n}\}$ converges to a domain $B$ (its kernel) bounded by a Jordan curve $C$. Let $\{f_{n}(\xi)\}$
denote a sequence of functions $f_{n}(\xi)$ such that, for each $n$, $f_{n}(0)=0$, $f_{n}^{'}(0)>0$ and $f_{n}(\xi)$
maps the unit disk $|\xi|<1$ onto the domain $B_{n}$. For the sequence $\{f_{n}(\xi)\}$ to converge uniformly on the
closed disk $|\xi|\le 1$ to a function $z=f(\xi)$ that vanishes at $0$, has positive first derivative at $0$, and
maps the open disk $|\xi|<1$ onto the domain $B$, it is necessary and sufficient that for every $\epsilon>0$, there
exists a number $N>0$ such that, for $n>N$, there exists a continuous one-to-one correspondence between the points
of the curves $C_{n}$ and $C$ such that the distance between any point of $C_{n}$ and the corresponding point of
$C$ will be less than $\epsilon$.

\begin{remark}
In the case of domains with arbitrary boundaries, not in particular of Jordan type, the question of convergence
in the closed disk, has been thoroughly investigated by Marku\v sevi\v c, \cite{m}.
\end{remark}
\noindent
We are now ready to tackle Problem number 2 that appear on the list of problems in section 5.

\section{The product of infinitely many shrinking factors}

Let $f(z)\in S(2,\alpha)$, and let $\{\phi_{n}\}_{n=1}^{\infty}$ be any sequence of real numbers. We defined
recursively: $\overline{g}_{\phi_{1}}$ is the Steiner symmetrization of $e^{i\phi_{1}}f(z)$, and 
$g_{\phi_{1}}=c_{\phi_{1}}\overline{g}_{\phi_{1}}$ where $c_{\phi_{1}}$ is the shrinking factor given by,
$$
c_{\phi_{1}}=\max\left\{\frac{1}{\overline{g}_{\phi_{1}}^{'}(0)},\frac{\alpha}{||\overline{g}_{\phi_{1}}||_{2}}\right\}.
$$
Next, $(\overline{g_{\phi_{1}}})_{\phi_{2}}$ is the Steiner symmetrization of $e^{i\phi_{2}}g_{\phi_{1}}(z)$, and
$(g_{\phi_{1}})_{\phi_{2}}=c_{\phi_{2}}(\overline{g_{\phi_{1}}})_{\phi_{2}}$, where $c_{\phi_{2}}$ is the shrinking
factor given by,
$$
c_{\phi_{2}}=\max\left\{\frac{1}{(\overline{g_{\phi_{1}}})_{\phi_{2}}^{'}(0)},
\frac{\alpha}{||(\overline{g_{\phi_{1}}})_{\phi_{2}}||_{2}}\right\}.
$$
The process proceeds indefinitely.

\begin{theorem}
Let $f(z)\in S(2,\alpha)$, then for any sequence $\{\phi_{n}\}_{n=1}^{\infty}$ of real numbers the limit function
$F=\lim_{n\rightarrow\infty}(\ldots((g_{\phi_{1}})_{\phi_{2}})\ldots)_{\phi_{n}}$ exists and the convergence is
uniform on compact subsets of $U$. The image $F(U)$ is a Steiner symmetric domain that includes $z=0$ and $F\in S(2,\alpha)$.
\end{theorem}
\noindent
{\bf Proof.} \\
Let us define $B_{n}=(\ldots((g_{\phi_{1}})_{\phi_{2}})\ldots)_{\phi_{n}}(U)$, $n=1,2,3,\ldots$. Then by the
definition of the recursive process, since $(\overline{\ldots((g_{\phi_{1}})_{\phi_{2}})\ldots})_{\phi_{n}}$ is
the Steiner symmetrization of $e^{i\phi_{n}}(\ldots((g_{\phi_{1}})_{\phi_{2}})\ldots)_{\phi_{n-1}}$, it follows
that the domain $\overline{B}_{n}=(\overline{\ldots((g_{\phi_{1}})_{\phi_{2}})\ldots})_{\phi_{n}}(U)$ is Steiner
symmetric, but $B_{n}=C_{\phi_{n}}\overline{B}_{n}$, a multiple by a number $0<c_{\phi_{n}}\le 1$ of a Steiner
symmetric domain. Hence $B_{n}$ itself is Steiner symmetric, and in particular the sequence $\{B_{n}\}_{n=1}^{\infty}$
is a sequence of simply connected domains, each of which contains the point $z=0$, and by the definition of
the Steiner symmetrization all of the domains $B_{n}$ contain the disk $|z|<{\rm dist}(0,\partial f(U))$. Hence,
the kernel of $\{B_{n}\}_{n=1}^{\infty}$, say $B$, exists and $\partial B$ is the faithful limit of $\partial B_{n}$ 
as $n\rightarrow\infty$. Hence $B_{n}\rightarrow B$ and by Theorem A (Carath\'eodory) the limit
$F=\lim_{n\rightarrow\infty}(\ldots((g_{\phi_{1}})_{\phi_{2}})\ldots )_{\phi_{n}}$ exists and is uniform on compact
subsets of $U$. Thus $B=F(U)$, and $F:\,U\rightarrow F(U)$ is conformal and satisfies the following normalization
$F(0)=0$, $F^{'}(0)>1$ and $||F||_{2}\ge\alpha$. This proves that $F\in S(2,\alpha)$ which is consistent
with the fact that $S(p,\alpha)$ is a compact family (recall the proof of Proposition 1.3). Moreover, the limit
$F(U)$ of the Steiner symmetric domains $\{(\ldots((g_{\phi_{1}})_{\phi_{2}})\ldots )_{\phi_{n}}(U)\}_{n=1}^{\infty}$
is Steiner symmetric. This is consistent with fact that $F(U)$ is a simply connected domain being the conformal
image of $U$. $\qed $ \\
\\
We can now give a sharp lower bound for the infinite product $\prod_{n=1}^{\infty} c_{\phi_{n}}$ of the shrinking
factors, that appears in the Problem 2 on the list of problems in section 5.

\begin{theorem}
1. Let $f\in S(2,\alpha)$ and let $\{\phi_{n}\}_{n=1}^{\infty}$ be any sequence of real numbers. Let 
$\{c_{\phi_{n}}\}_{n=1}^{\infty}$ be the corresponding sequence of the shrinking factors. Then we
have the following estimate:
$$
\max\left\{1,\alpha\right\}\cdot\left\{\frac{1}{\pi}\int_{0}^{2\pi}\Re\{f(e^{i\theta}\}
\Re\{e^{i\theta}f^{'}(e^{i\theta})\}d\theta\right\}^{-1/2}\le\prod_{n=1}^{\infty}c_{\phi_{n}}\le 1,
$$
and these bounds on $\prod_{n=1}^{\infty}c_{\phi_{n}}$ are sharp bounds. \\
2. Let us define recursively the following sequence of mappings: $f_{\phi_{1}}$ is the Steiner symmetrization
of $e^{i\phi_{1}}f$. For $n\in\mathbb{Z}^{+}$, let $(\ldots((f_{\phi_{1}})_{\phi_{2}})\ldots)_{\phi_{n+1}}$ be
the Steiner symmetrization of $e^{i\phi_{n+1}}(\ldots((f_{\phi_{1}})_{\phi_{2}})\ldots)_{\phi_{n}}$. Then the
limit function $G(z)=\lim_{n\rightarrow\infty}(\ldots((f_{\phi_{1}})_{\phi_{2}})\ldots)_{\phi_{n}}(z)$ exists
and is uniform on compact subsets of $U$. Moreover, we have the following identity:
$$
\prod_{n=1}^{\infty}c_{\phi_{n}}=\max\left\{\frac{1}{G^{'}(0)},\frac{\alpha}{||G||_{2}}\right\}.
$$
\end{theorem}

\begin{remark}
Theorem 9.2 gives some kind of solution to Problem 2 on the list of problems that appear in
section 5.
\end{remark}
\noindent
{\bf A proof of Theorem 9.2.} \\
1. We will use the recursive sequence $\{f_{\phi_{n}}\}_{n=1}^{\infty}$, that was defined in part 2 of
Theorem 9.2. Then as in the proof of Theorem 9.1 that dealt with the sequence 
$\{(\ldots((g_{\phi_{1}})_{\phi_{2}})\ldots)_{\phi_{n}}\}_{n=1}^{\infty}$,
based on Theorem A (Carath\'eodory), the limit $G=\lim_{n\rightarrow\infty}(\ldots((f_{\phi_{1}})_{\phi_{2}})\ldots)_{\phi_{n}}$
exists and is uniform on compact subsets of $U$. We note that the recursive process outlined by the newer sequence
$\{(\ldots((f_{\phi_{1}})_{\phi_{2}})\ldots)_{\phi_{n}}\}_{n=1}^{\infty}$ is simpler than the original recursive
process that was described on section 3, in that we do not multiply by the shrinking factors $c_{\phi_{n}}$ after each symmetrization
was done. The purpose of those multiplications was to optimize, i.e. make as small as possible, each element of the
sequence of functions produced. That without leaving the family $S(2,\alpha)$. We will soon see that if our goal was to
optimize the limiting function and not the each element of the sequence, then this can be accomplished by a single
multiplication by just one shrinking factor. In fact this is the key idea for the current proof.

The first step is to note that each element of the old sequence, the $g$-sequence, is a multiple by a constant of the
corresponding element of the new system, the $f$-sequence. The constant, though is not a single shrinking factor. Clearly
by the definitions of $g_{\phi_{1}}$ and of $f_{\phi_{1}}$ we have the formula $g_{\phi_{1}}=c_{\phi_{1}}f_{\phi_{1}}$.
Next, $(f_{\phi_{1}})_{\phi_{2}}$ is the conformal mapping $U\rightarrow (e^{i\phi_{2}}f_{\phi_{1}}(U))^{*}$ where
$A^{*}$ denotes the Steiner symmetrization of the domain $A$. This conformal mapping is normalized as follows:
$(f_{\phi_{1}})_{\phi_{2}}(0)=0$ and $(f_{\phi_{1}})_{\phi_{2}}^{'}(0)>0$. That definition of $(f_{\phi_{1}})_{\phi_{2}}$
should be compared with the definition of $(g_{\phi_{1}})_{\phi_{2}}$ which equals the shrinking factor, $c_{\phi_{2}}$
multiplying the conformal mapping $U\rightarrow (e^{i\phi_{2}}g_{\phi_{1}}(U))^{*}$. But we already have the formula
$e^{i\phi_{2}}g_{\phi_{1}}(U)=e^{i\phi_{2}}(c_{\phi_{1}}f_{\phi_{1}}(U))=c_{\phi_{1}}(e^{i\phi_{2}}f_{\phi_{1}}(U))$. Thus
we deduce that the Steiner symmetrization of $e^{i\phi_{2}}g_{\phi_{1}}(U)$ equals to $c_{\phi_{1}}$ times the Steiner
symmetrization of $e^{i\phi_{2}}f_{\phi_{1}}(U)$. In other words the relation between the image of the first conformal
mapping to the image of the second conformal mapping is multiplication by the shrinking factor $c_{\phi_{1}}$, where
we recall that $0<c_{\phi_{1}}\le 1$. Hence by composition of conformal mappings we get the following (second) formula
$(g_{\phi_{1}})_{\phi_{2}}=c_{\phi_{2}}c_{\phi_{1}}(f_{\phi_{1}})_{\phi_{2}}$. Similarly, the general case follows by
an inductive argument. We obtain the general formula:
$$
(\ldots((g_{\phi_{1}})_{\phi_{2}})\ldots)_{\phi_{n}}=\left(\prod_{k=1}^{n}c_{\phi_{k}}\right)\cdot
(\ldots((f_{\phi_{1}})_{\phi_{2}})\ldots)_{\phi_{n}},\,\,\,\,\,\forall\,n\in\mathbb{Z}^{+}.
$$
Passing to the limit $n\rightarrow\infty$ gives us:
$$
\lim_{n\rightarrow\infty}(\ldots((g_{\phi_{1}})_{\phi_{2}})\ldots)_{\phi_{n}}=\left(\prod_{k=1}^{\infty}c_{\phi_{k}}\right)\cdot
(\ldots((f_{\phi_{1}})_{\phi_{2}})\ldots)_{\phi_{n}},
$$
or simply (using our notations for the limits):
\begin{equation}
\label{eq11}
F=\left(\prod_{k=1}^{\infty}c_{\phi_{k}}\right)\cdot G.
\end{equation}
If the sequence $\{\phi_{n}\}_{n=1}^{\infty}$ deforms $f(U)$ to a disk $D(0,R)$ in the newer process, then this disk
has an area which equals the area of $f(U)$. We obtain the following equation with the unknown $R$:
$$
\pi R^{2}=\int_{0}^{2\pi}\Re\{f(e^{i\theta})\}\Re\{e^{i\theta}f^{'}(e^{i\theta})\}d\theta.
$$
Hence:
$$
R=\left\{\frac{1}{\pi}\int_{0}^{2\pi}\Re\{f(e^{i\theta})\}\Re\{e^{i\theta}f^{'}(e^{i\theta})\}\right\}^{1/2}.
$$
This means that for this particular sequence $\{\phi_{n}\}_{n=1}^{\infty}$, We obtain the very simple formula for
the conformal mapping $G$, namely that $G:\,U\rightarrow D(0,R)$, given by $G(z)=R\cdot z$. In particular we get
$G^{'}(0)=||G||_{2}=R$ and hence the corresponding shrinking factor, which optimize $G$ is given by:
$$
\max\left\{\frac{1}{R},\frac{\alpha}{R}\right\}=\max\left\{1,\alpha\right\}\cdot\frac{1}{R}=
\max\left\{1,\alpha\right\}\left\{\frac{1}{\pi}\int_{0}^{2\pi}\Re\{f(e^{i\theta})\}\Re\{e^{i\theta}f^{'}(e^{i\theta})\}\right\}^{-1/2}.
$$
This concludes the proof of the inequality of part 1, including its sharpness. \\
2. In the general case the limiting function for an arbitrary sequence of real numbers $\{\phi_{n}\}_{n=1}^{\infty}$ is
the conformal mapping $G:\,U\rightarrow G(U)$, where $G(0)=0$ and $G^{'}(0)>0$. Now the general shrinking factor (of $G$) is
given by:
$$
c=\max\left\{\frac{1}{G^{'}(0)},\frac{\alpha}{||G||_{2}}\right\}.
$$
On the other hand, by equation \ref{eq11} this shrinking factor is given by the infinite product:
$$
c=\prod_{n=1}^{\infty}c_{\phi_{n}}.
$$
This concludes the proof of part 2 of our theorem. $\qed $ \\

\begin{corollary}
Let $f\in S(2,\alpha)$ and let $\{\phi_{n}\}_{n=1}^{\infty}$ be any sequence of real numbers. Let $G$ be the
limiting function of the newer recursive process, i.e. $G=\lim_{n\rightarrow\infty}(\ldots((f_{\phi_{1}})_{\phi_{2}})\ldots)_{\phi_{n}}$.
Then we have the sharp estimate:
$$
\max\left\{1,\alpha\right\}\cdot\left\{\frac{1}{\pi}\int_{0}^{2\pi}\Re\{f(e^{i\theta})\}
\Re\{e^{i\theta}f^{'}(e^{i\theta})\}d\theta\right\}^{-1/2}\le\max\left\{\frac{1}{G^{'}(0)},\frac{\alpha}{||G||_{2}}\right\}\le 1.
$$
In particular, if $\alpha$ is taken to be small enough, then:
$$
G^{'}(0)\le\left\{\frac{1}{\pi}\int_{0}^{2\pi}\Re\{f(e^{i\theta})\}\Re\{e^{i\theta}f^{'}(e^{i\theta})\}d\theta\right\}^{1/2}.
$$
This last upper bound is sharp.
\end{corollary}

\noindent
{\it Ronen Peretz \\
Department of Mathematics \\ Ben Gurion University of the Negev \\
Beer-Sheva , 84105 \\ Israel \\ E-mail: ronenp@math.bgu.ac.il} \\ 
 
\end{document}